%% file: main.tex
\documentclass[12pt]{article}

\usepackage[english]{babel}
\usepackage[utf8]{inputenc}
\usepackage[T1]{fontenc}
\usepackage{lmodern}
\usepackage[a4paper, margin=2cm]{geometry}
\usepackage[onehalfspacing]{setspace}
\usepackage{microtype}
\usepackage{enumitem}
\setlist[itemize]{label=\textbullet}

\usepackage{graphicx}
\graphicspath{{figures/}}
\usepackage{booktabs}
\usepackage{tabularx}
\usepackage{multirow}
\usepackage{array}
\usepackage{caption}
\usepackage{float}
\usepackage{xcolor}

\usepackage{mathtools}
\usepackage{amssymb}
\usepackage{amsthm}

\newcommand{\sgap}{\ensuremath{\varrho}}   
\newcommand{\ascore}{\ensuremath{\Lambda}} 
\newcommand{\tlim}{\ensuremath{\Gamma}}    

\usepackage{csquotes}
\usepackage[
  backend=biber,
  style=alphabetic,
  sorting=none,
  language=american,
  maxcitenames=1,
  mincitenames=1,
  maxbibnames=7,
  minbibnames=7
]{biblatex}
\newcommand{\keywords}[1]{\textbf{Keywords}: \itshape #1 \normalfont}

\title{\Large Anytime Solver Evaluation with a Normalized Signed Primal Integral and Explicit Reference Policies --- Extended Version}
\author{
  \normalsize Florian Rascoussier\textsuperscript{1,2}\thanks{ORCID identifiers: \begin{tabular}[t]{@{}l@{\quad}l@{}}Florian Rascoussier: \href{https://orcid.org/0009-0005-3253-9814}{0009-0005-3253-9814} & Romain Billot: \href{https://orcid.org/0000-0002-2614-3822}{0000-0002-2614-3822}\\ Lina Fahed: \href{https://orcid.org/0000-0002-7207-0580}{0000-0002-7207-0580} & Christine Solnon: \href{https://orcid.org/0000-0002-0919-496X}{0000-0002-0919-496X}\end{tabular}}\quad
  Romain Billot\textsuperscript{1}\quad
  Lina Fahed\textsuperscript{1}\quad
  Christine Solnon\textsuperscript{2}\\[0.5ex]
  \small \textsuperscript{1}IMT Atlantique, Lab-STICC, CNRS, UMR 6285 (équipe DECIDE), France\\
  \small \textsuperscript{2}INSA Lyon, Inria, CITI, UR3720, 69621 Villeurbanne, France (équipe EMERAUDE)
}
\date{Version 3, 19/10/2026}

\usepackage[
    colorlinks=true,
    linkcolor=blue,
    citecolor=blue,
    urlcolor=blue
]{hyperref}

\hypersetup{pdfauthor={Florian Rascoussier, Romain Billot, Lina Fahed, Christine Solnon}}

\begin{document}

\maketitle

\input{sections/0-core-content}

\printbibliography

\clearpage
\appendix
\input{sections/1-appendices}

\end{document}

%% file: sections/0-core-content.tex

\vspace{-0.5cm} 

\begin{abstract}
Anytime solvers return usable solutions before they terminate and improve them while time remains. Their progress is commonly summarized by a primal integral, a final gap, or convergence curves. Each summary leaves consequential choices open: how to score a run before its first feasible solution, whether poor incumbents are truncated, and whether the reference value is updated after the experiment or version-frozen beforehand. These choices become visible when runs produce no valid incumbent or improve a published best-known value.

We study a normalized signed primal integral built from a bounded relative gap. It assigns an intrinsic worst value to an empty run, requires no acceptance threshold, and assigns negative instantaneous gaps to incumbents that beat a frozen reference. We compare the smooth difference-over-sum kernel with a signed version of Berthold's max-normalized gap and use the former as a working default. We also distinguish analysis-time from version-frozen reference policies and recommend reading the score with the raw final gap, mean convergence curve, and target-attainment curve.

We evaluate these choices on a five-arm routing campaign and two model-fidelity ladders. The two signed kernels preserve every panel ordering in this study, whereas a common acceptance threshold compresses the distances between arms and reverses one panel ordering. Updating 54 of the 212 reference values changes score levels and removes all negative scores, while leaving the observed panel orderings unchanged. An exploratory screen also identifies 12 panel comparisons in which similar integral scores conceal materially different endpoints or attainment rates. The implementation, frozen inputs, and generators are openly released.

\vspace{2ex}

\keywords{Anytime optimization, Performance evaluation, Primal integral, Reference policy, Benchmarking}
\end{abstract}

\setcounter{tocdepth}{2}
\tableofcontents
\clearpage

\section{Introduction}\label{sec:intro}

An optimization procedure is \emph{anytime} when it can be interrupted and still return a usable answer whose quality tends to improve with additional computation \parencite{zilbersteinUsingAnytimeAlgorithms1996}. Branch-and-bound procedures, construction heuristics, local searches, and modern hybrid solvers all expose this behavior, although their trajectories differ sharply. One method may spend most of its budget without a feasible solution and then improve quickly; another may construct a good solution immediately and plateau. A final-value table records only the last point of each trajectory and cannot distinguish these behaviors.

The primal integral introduced by \textcite{bertholdMeasuringImpactPrimal2013} is the established scalar response to that problem. It integrates the best-so-far gap over time and is used in vehicle-routing evaluation through the 12th DIMACS Implementation Challenge and later benchmarking environments \parencite{DIMACS2021_rules,thyssensRoutingArenaBenchmark2023,bertoRL4COExtensiveReinforcement2025}. Mean convergence curves and target-attainment views provide complementary readings of the same progress. The need for an anytime instrument is therefore not the gap addressed here. The open issue is its \emph{instantiation}: a score or curve still requires a gap function, a reference policy, a value before the first incumbent, a convention for invalid or empty runs, and an aggregation rule.

These choices matter precisely in the difficult cases for which anytime evaluation is most useful. A raw relative gap is unbounded before the first solution. Competition scores often impose an acceptance threshold, making a finite but poor incumbent indistinguishable from no solution. A reference chosen from a published best-known-solution store may be beaten during the experiment, so a measure designed for a post-hoc reference must either be extended below zero or applied outside its original setting. Finally, a scalar deliberately compresses a trajectory and cannot alone state where a run ended or how frequently a target was attained.

This paper studies one coherent response to those cases. The instantaneous squeezed gap $\sgap$ is finite and bounded for every positive minimization objective, with the no-incumbent state fixed at its intrinsic supremum. Under a frozen reference, incumbents that improve the reference make negative instantaneous contributions to its budget-normalized time integral $\ascore$. Crucially, the reference itself is treated as a methodological choice between two regimes: an analysis-time reference computed after every run is known, and a versioned snapshot frozen before the experiment. The raw final gap, mean convergence curve, and goal-attainment curve are recommended companions because they recover information the scalar compresses.

Three contributions structure the paper. First, we name and analyze the \emph{Normalized Signed Primal Integral}, state its domain and exact finite-sum evaluation, and compare its difference-over-sum kernel with the signed max-form variant that exactly matches Berthold above the reference. We adopt the former as a smooth practical default while explicitly recognizing that the latter is a useful alternative. Second, we make the reference policy explicit and show how the same trajectories acquire different signs and score levels under analysis-time and frozen references. The key property of this second regime is making the selected per-instance reference set precisely versionable, with later improvements producing a new snapshot rather than silently altering prior scores. Each value can be pinned to a named release or commit of an evolving BKS store such as MAMUT-routing, or to a public benchmark repository such as the GitLab reference repository accompanying \textcite{blauthVehicleRoutingTimedependent2024}.\footnote{MAMUT-routing version history: \url{https://mamut-routing.onyr.net/history/index.html}; historical VRP BKS records at CombOpt: \url{http://combopt.org/history/}; Blauth et al.'s public reference repository: \url{https://gitlab.com/muelleratorunibonnde/vrptdt-benchmark}.}  Third, we evaluate the practical consequences on 8,800 independently validated runs and two model-fidelity ladders, quantifying failure conventions, metric sensitivity, reference sensitivity, and the information recovered by companion views. Every result regenerates deterministically from machine-checkable archived experiment artifacts.

The definitions apply to any anytime procedure minimizing a strictly positive objective. Vehicle routing supplies the empirical setting because it combines stochastic solvers, evolving best-known-solution (BKS) stores, difficult feasibility, and independently checkable incumbents. Section~\ref{sec:background} establishes the design problem and prior art. Section~\ref{sec:metric} defines the metric and its alternative. Section~\ref{sec:policies} develops the two reference regimes and companion views. Section~\ref{sec:evaluation} answers three empirical research questions. Section~\ref{sec:practice} gives the minimal reporting, aggregation, limitation, and reproducibility guidance before Section~\ref{sec:conclusion} concludes.

\section{Background and Design Problem}\label{sec:background}\label{sec:related}

\subsection{Runs, trajectories, and two empirical views}

Considering an experiment campaign, a randomized solver configuration induces a distribution over runs. The object observed in one run is a timestamped incumbent trajectory, while the object compared across seeds or instances is a statistic of a collection of trajectories. The distinction matters because averaging candidate solutions, averaging gap curves, and averaging per-run areas need not express the same estimand unless the operations are chosen to commute.

Two classical empirical views answer different questions \parencite{hoosStochasticLocalSearch2004,hansenCOCOPlatformComparing2021,lopez-ibanezUsingEmpiricalAttainment2025}. A fixed-budget convergence curve asks how solution quality evolves as time is increased. A fixed-target attainment curve asks what fraction of runs has reached a declared quality goal by each time. For an individual run, an ordinary incumbent-quality curve has no finite value before the first feasible solution, unless a convention is specified for that state. An attainment curve remains defined when a run never reaches the goal: its individual hitting time is right-censored or infinite, and it contributes zero attainment at every finite time.

The two views are complementary to a scalar area measure. An integral compares whole trajectories at one declared horizon; an attainment curve keeps the distributional meaning of reaching a target; a final gap answers the endpoint question directly. Our methodological contribution is to use compatible gap and failure conventions so that a run does not disappear from one view while receiving an arbitrary penalty in another.

\subsection{Berthold and DIMACS}

For an incumbent objective value $z$ and an optimal or best-known reference $z^*$, Berthold's primal gap on a strictly positive minimization objective reduces to
\begin{equation}
  \gamma_B(z,z^*) = \frac{z-z^*}{z}=1-\frac{z^*}{z}, \qquad z\ge z^*,
  \label{eq:berthold-positive-gap}
\end{equation}
with gap $1$ before any incumbent exists \parencite{bertholdMeasuringImpactPrimal2013}. His general definition is the absolute form $\gamma_B(z,z^*)=\lvert z-z^*\rvert/\max(z,z^*)$, which the reduction above instantiates whenever the reference is not beaten. The primal integral is the area under the resulting best-so-far gap curve. Its reference is normally the best value known at analysis time and is updated if a run improves it, preserving $z\ge z^*$ and a nonnegative gap. The confined primal integral changes the time weighting to emphasize different parts of the run, an orthogonal design choice considered in Appendix~\ref{app:logtime} \parencite{bertholdConfinedPrimalIntegral2021}.

The DIMACS vehicle-routing score is direct prior art for a signed, budget-normalized integral against a frozen reference \parencite{DIMACS2021_rules}. It initializes the trajectory at $\mathit{ACC}\,z^*$ for an acceptance threshold $\mathit{ACC}=1.1$, accepts only solutions below $\mathit{ACC}\,z^*$, integrates the raw excess gap $g = z/z^*-1$, and multiplies the average by 100. A run that never crosses the threshold scores 10 whether it found nothing or improved steadily above $1.1z^*$. A solution below the frozen best-known value produces a negative contribution. Routing Arena and RL4CO subsequently adopted the same rule \parencite{thyssensRoutingArenaBenchmark2023,bertoRL4COExtensiveReinforcement2025}.

Our contribution is to isolate the instantaneous normalization, acceptance threshold, pre-incumbent value, reference policy, and curve conventions; analyze a threshold-free signed default and measure when these choices affect conclusions.

\subsection{Neighboring benchmarking constructions}

Performance profiles aggregate ratios to the best solver time over instances, while data profiles aggregate the fraction of problems solved within normalized effort \parencite{dolanBenchmarkingOptimizationSoftware2002,moreBenchmarkingDerivativeFree2009}. COCO/BBOB expected running time and target ECDFs are nonnegative fixed-target constructions over known target grids \parencite{hansenCOCOPlatformComparing2021}. IOHanalyzer's area over the convergence curve and empirical attainment functions establish close links between area and attainment views \parencite{wangIOHanalyzerDetailedPerformance2022,lopez-ibanezUsingEmpiricalAttainment2025}. Hypervolume-over-time formulations instead score a time-quality front against a reference point \parencite{lopez-ibanezAutomaticallyImprovingAnytime2014}. These methods answer related benchmarking questions, but none removes the need to define what a missing incumbent or a beaten external reference means in an area-under-gap score.

The primal-dual integral is also distinct. It uses a solver's own dual bound and is therefore self-referential when solvers have different bounding strength. The present setting compares heuristic and exact procedures through a shared external objective evaluator and a shared per-instance reference. The external reference is what makes run-level values commensurable, and its provenance becomes part of the experiment.

\subsection{The timing and validation contract}

All studies below measure elapsed wall-clock time. This choice matches the operational question addressed by anytime evaluation: what solution quality is available when a real deadline expires? It also provides a common basis for comparing solvers that use different numbers of threads, for which CPU time would account for parallel work differently. The DIMACS rules adopt the same convention \parencite{DIMACS2021_rules}. Wall-clock time is not universally preferable, however. CPU time may provide the more informative measure in a strictly sequential study of a single solver, and the metric itself accommodates either clock.

A clock becomes meaningful only as part of an explicit timing protocol. A study should therefore state the horizon, the work included within it, the timestamp granularity, the thread allocation, and the machine class. In particular, it should clarify whether initialization, model construction, and validation are charged to the budget. In our experiments, the horizon is announced to each solver in advance, so a solver may legitimately adapt its search to the available time. The resulting score describes performance conditional on that declared horizon, rather than a budget-independent property of the solver.

Solver-reported feasibility and objective values do not directly define the measured trajectory. Each emitted candidate is independently checked against the original problem, and its objective is recomputed before the candidate can enter the best-so-far envelope at its recorded timestamp. This distinction is essential when a solver works with an approximation of the evaluated problem. For example, the Hexaly Tier-2 slice ladders used in Section~\ref{sec:evaluation} contain solutions that are feasible for the solver's discretized model but infeasible for the original time-dependent functions. Such candidates create no improvement event. If a run never produces a checker-valid candidate, its trajectory remains at the pre-incumbent value for the full horizon. We retain these runs in every aggregate and report the numbers of rejected candidates and empty streams separately. The bounded score gives failures a well-defined contribution, but the separate counts remain necessary to distinguish a failure to produce valid output from poor performance on valid solutions.

\subsection{Design choices and their consequences}

Table~\ref{tab:design-choices} separates choices that are sometimes hidden inside the phrase \enquote{primal integral}. They need not be resolved identically in every study. They must, however, be resolved consistently across the score, curves, and aggregation if those outputs are to describe the same experiment.

\begin{table}[htbp]
  \centering
  \small
  \caption{Design choices in an anytime area measure. The rightmost column states the convention evaluated in this paper, not a universal requirement.}
  \label{tab:design-choices}
  \begin{tabular}{@{}>{\raggedright\arraybackslash}p{0.21\linewidth}>{\raggedright\arraybackslash}p{0.31\linewidth}>{\raggedright\arraybackslash}p{0.40\linewidth}@{}}
    \toprule
    Choice & Alternatives & Convention studied here \\
    \midrule
    Instantaneous gap & Raw, max-normalized, difference-over-sum, or thresholded & Compare the alternatives; use the bounded signed difference-over-sum kernel as the working default \\
    Before feasibility & Undefined, imputed, initialized by a finite construction, or intrinsic worst value & Assign the kernel supremum 1 until the first checker-valid incumbent \\
    Reference timing & Best value at analysis time or a frozen earlier snapshot & Treat both as first-class policies and never pool scores across snapshots \\
    Invalid candidates & Trust solver status, penalize, discard the run, or validate externally & Exclude invalid candidates from the envelope; retain and count empty valid streams \\
    Time and horizon & CPU or wall clock; one or several budgets & Declare the clock and horizon; normalize area by the chosen horizon \\
    Across-run summary & Mean, median, quantiles, or robust alternatives & Use the arithmetic mean for an expected-performance estimand and its area identity \\
    Information beyond area & Endpoint, convergence, and target-attainment summaries & Recommend the views needed to interpret the claims made by the study \\
    \bottomrule
  \end{tabular}
\end{table}

The objective domain is another substantive choice. A bounded relative kernel does not remove the dependence on an objective representation: multiplicative rescaling is harmless, while an additive shift is not. We therefore restrict the formal development to strictly positive minimization objectives. Similarly, the horizon is normalized but not eliminated. A score at five minutes and a score at one hour answer different interruption questions even though both lie on the same numerical scale.

\section{The Normalized Signed Primal Integral}\label{sec:metric}

\subsection{Incumbent trajectory and squeezed gap}\label{subsec:metric-defs}

Consider a run over a declared horizon $[0,\tlim]$. After validation and best-so-far filtering, its incumbent events are $(t_1,z_1),\ldots,(t_m,z_m)$ with $0\le t_1\le\cdots\le t_m\le\tlim$ and $z_1>\cdots>z_m>0$. The right-continuous value $z(t)$ equals the current incumbent between events and $+\infty$ before the first valid incumbent. If execution stops before $\tlim$, including after a crash, the last checker-valid incumbent is held until $\tlim$; a stopped run with no valid event is an empty run and has $z(t)=+\infty$ throughout.

For a strictly positive reference $z^*$, define the squeezed gap
\begin{equation}
  \sgap(z,z^*)=\frac{z-z^*}{z+z^*}=\frac{g}{g+2}, \qquad g=\frac{z-z^*}{z^*},
  \label{eq:squeezed-gap}
\end{equation}
with $\sgap(+\infty,z^*)=1$. Writing $x=z/z^*$ gives
\begin{equation}
  \sgap(z,z^*)=\frac{x-1}{x+1}=\tanh\!\left(\frac{1}{2}\ln x\right).
  \label{eq:squeezed-tanh}
\end{equation}
The kernel is strictly increasing, scale-invariant, and log-symmetric: multiplying $z$ and $z^*$ by the same positive constant changes nothing, while reciprocal ratios produce opposite gaps. It is not translation-invariant, which is why positive objectives are a real scope condition rather than a convenience.

For admissible finite objectives $z>0$, the range is $(-1,1)$; reference equality gives zero, and $1$ is approached only as the incumbent diverges. The extended no-incumbent state is assigned $1$, yielding the combined range $(-1,1]$. If $z=0$ were admitted, $-1$ would be attainable and the range would be $[-1,1]$. Objectives that may be zero or negative require a declared positive transformation, whose influence must be analyzed because additive shifts change the metric.

\begin{figure}[htbp]
  \centering
  \includegraphics[width=\linewidth]{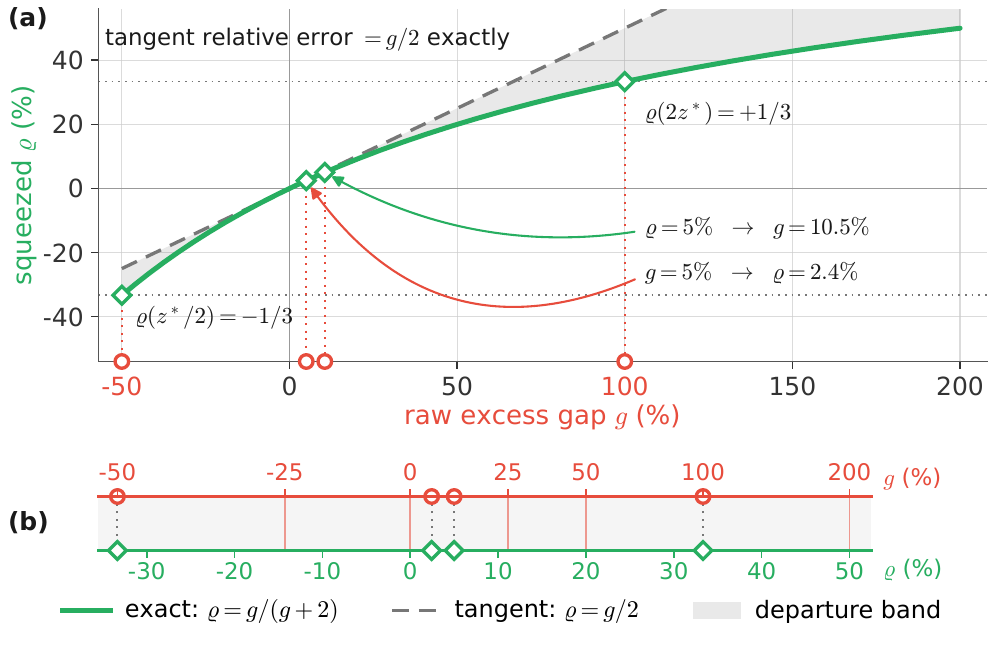}
  \caption{Conversion between the raw excess gap $g$ and the squeezed gap $\sgap=g/(g+2)$. (a) The exact conversion against its tangent $\sgap=g/2$ at the reference; the shaded band is the departure between the two, and the diamonds mark the landmarks $\sgap(2z^*)=+1/3$ and $\sgap(z^*/2)=-1/3$. (b) The two scales as one dual ruler. The squeezed scale compresses large raw gaps while remaining linear to first order near the reference.}
  \label{fig:scale-conversion}
\end{figure}

Near the reference, $\sgap=g/(g+2)\approx g/2$. Thus a raw 5\% excess gap maps to approximately 2.4\% on the squeezed scale, while a squeezed 5\% corresponds to a raw gap of approximately 10.5\%. Figure~\ref{fig:scale-conversion} draws the exact conversion, its tangent at the reference, and the two scales side by side; Appendix~\ref{app:math} tabulates landmark values. Labels and goal definitions must distinguish the two scales. Goal levels in this paper are stated on the familiar raw scale and converted internally by $\sgap=g/(g+2)$.

The difference-over-sum kernel has precedents in symmetric percentage errors, normalized-difference indices, and Michelson contrast \parencite{hyndmanAnotherLookMeasures2006,tuckerRedPhotographicInfrared1979,michelsonStudiesOptics1927}. Neither the algebraic form nor signed area integration is claimed as novel. Its value here is operational: one analytic expression covers both sides of the reference and tends smoothly to an intrinsic worst value.

\subsection{The score and exact evaluation}

The \emph{Normalized Signed Primal Integral} is the time average of the squeezed gap:
\begin{equation}
  \ascore(\tlim)=\frac{1}{\tlim}\int_0^{\tlim}\sgap\!\left(z(t),z^*\right)\,\mathrm{d}t
  =\frac{1}{\tlim}\sum_{i=0}^{m}\sgap(z_i,z^*)\,(t_{i+1}-t_i),
  \label{eq:ascore}
\end{equation}
where $t_0=0$, $z_0=+\infty$, and $t_{m+1}=\tlim$. On a timestamped incumbent trace this is an exact finite left sum, not an interpolation. Equal timestamps are coalesced by keeping the best objective, events beyond the horizon are discarded, and invalid candidates never enter the envelope.

Because every instantaneous value lies in $(-1,1]$, the score also lies in $(-1,1]$. An empty run scores exactly 1. A run reaching the reference immediately and never improving scores 0. Smaller is better. The horizon remains part of the estimand: changing $\tlim$ can change rankings because it asks a different operational question. Normalization removes the time unit but does not make the score budget-independent.

The sign requires careful interpretation. A negative $\ascore$ means that the magnitude-weighted signed area below the reference exceeds the area above it. It does not imply that the run spent most of its time below the reference, because gap magnitude and duration both contribute. Conversely, finishing below the reference does not imply $\ascore<0$: a run may cross the reference near the deadline after accumulating a large positive area. Discovery of a better endpoint is therefore reported through the final gap and trajectory, not inferred from the integral's sign.

\subsection{A valid alternative: the signed max-form}\label{subsec:maxform}

The direct signed extension of Berthold's positive-objective gap is
\begin{equation}
  \sgap_B(z,z^*)=\frac{z-z^*}{\max(z,z^*)}.
  \label{eq:maxform}
\end{equation}
It equals Berthold's gap exactly for $z\ge z^*$, equals the raw excess gap below the reference, is scale-invariant and log-symmetric, and uses the same no-incumbent value 1. It is $C^1$ at $z=z^*$; only its second derivative is discontinuous at the reference.

\begin{figure}[htbp]
  \centering
  \includegraphics[width=0.88\linewidth]{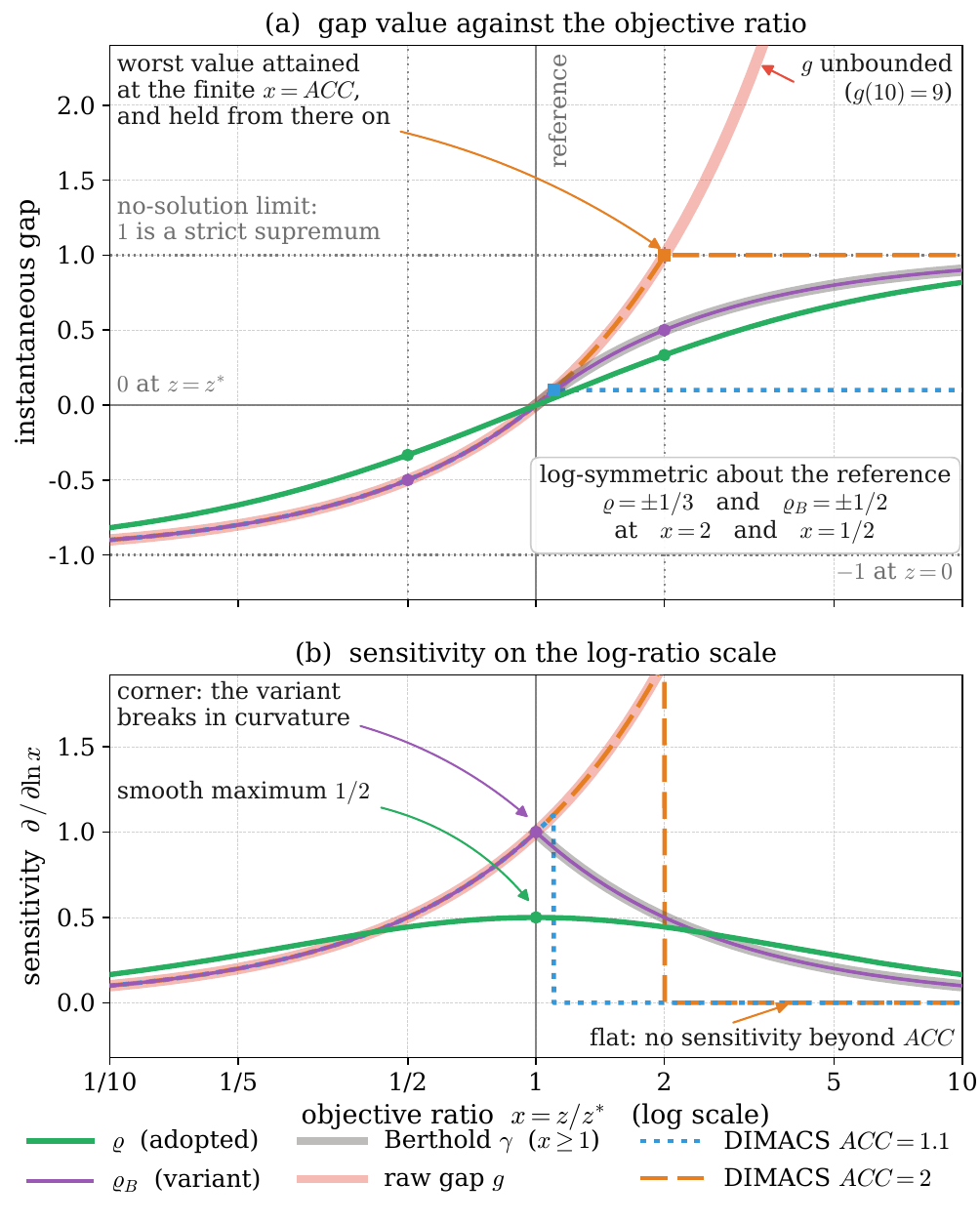}
  \caption{Instantaneous gap choices as functions of the objective ratio $z/z^*$. The squeezed gap $\sgap$ and signed max-form $\sgap_B$ are both bounded, signed, scale-invariant, and log-symmetric. The signed max-form coincides with Berthold above the reference and changes curvature at the reference; the squeezed gap uses one smooth difference-over-sum expression. DIMACS thresholds flatten finite regions and therefore make some improving trajectories indistinguishable from empty runs.}
  \label{fig:gap-shapes}
\end{figure}

We write $\ascore_B$ for the score of Equation~(\ref{eq:ascore}) computed with $\sgap_B$ in place of $\sgap$. The signed max-form has two genuine advantages. It preserves the familiar local scale $\sgap_B\approx g$, and in the classic nonnegative regime its integral is exactly Berthold's normalized primal integral. The squeezed gap instead uses one smooth analytic kernel across the whole positive domain and a single closed-form inverse. These are different design preferences, as illustrated in Figure~\ref{fig:gap-shapes}, not a proof that one candidate is correct and the other defective.

We use $\sgap$ as the working default because the frozen-reference regime is central here and a single smooth expression makes its two sides visually and analytically uniform. The alternative remains equally valid. The empirical comparison in Section~\ref{sec:evaluation} finds no Tier-1 panel-pair ordering reversal between the two signed kernels: on these data the choice changes scale and some run-level signs but not the solver ordering. Both kernels thus preserve a similar semantic treatment of empty and below-reference intervals and support similar conclusions. Appendix~\ref{app:maxform} records the derivative and range details.

\subsection{A compact worked example}\label{subsec:worked-example}

\begin{figure}[htbp]
  \centering
  \includegraphics[width=\linewidth]{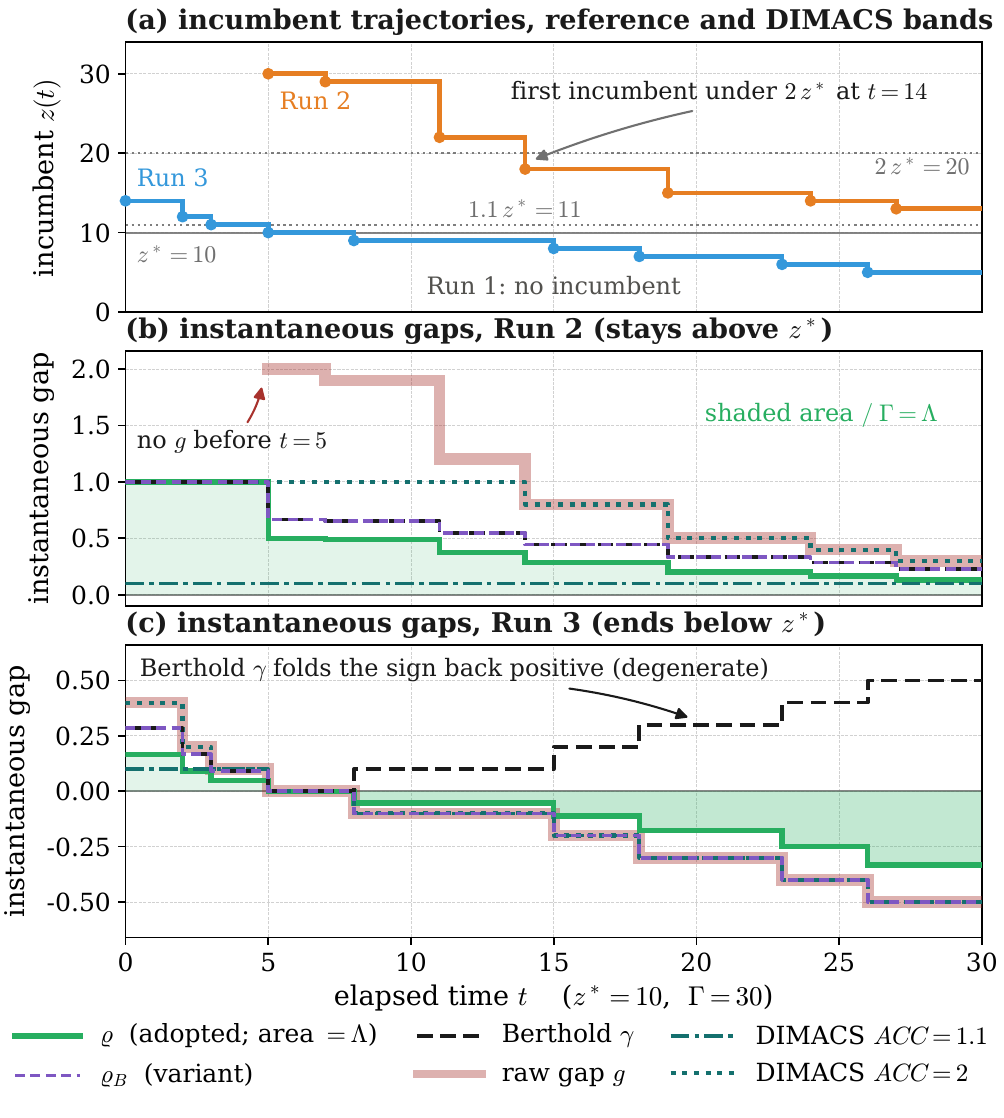}
  \caption{The worked example drawn. (a) The three runs against the reference $z^*=10$ and the two DIMACS acceptance levels: Run~1 produces no incumbent and Run~2 never enters the $1.1\,z^*$ band. (b, c) The instantaneous gaps of Runs 2 and 3 under the compared instantiations; the shaded signed area of $\sgap$, normalized by $\tlim$, is $\ascore$, with the below-reference part (darker) counting negative. Run~1 needs no gap panel: every bounded instantiation sits at its worst value throughout, and the raw gap has no finite value at all. Overlapping dash patterns mark exact coincidences: each instantiation holds its pre-incumbent value until its first accepted incumbent ($1$ for the bounded kernels, $\mathit{ACC}-1$ for DIMACS), $\sgap_B$ equals Berthold's gap above the reference and the raw gap below it, and each DIMACS integrand equals the raw gap once its threshold admits the incumbent. Berthold's gap below the reference in panel~(c) is behavior outside the measure's design envelope, shown deliberately.}
  \label{fig:worked-example}
\end{figure}

Three trajectories against $z^*=10$ over a horizon $\tlim=30$ make the conventions concrete. Run~1 produces no valid incumbent. Run~2 first becomes feasible at $t=5$, improves repeatedly from 30 to 13, but never reaches the reference. Run~3 is feasible at the origin, reaches the reference at $t=5$, and then improves to 5. Figure~\ref{fig:worked-example} draws the trajectories and their instantaneous gaps under the comparison set of Figure~\ref{fig:gap-shapes}; Table~\ref{tab:pi-example} evaluates the exact finite sums. Away from the reference the signed max-form assigns larger magnitudes, its normalizer $\max(z,z^*)$ being smaller than the sum $z+z^*$, but it gives the same ordering and the same treatment of empty and below-reference intervals as $\sgap$.

\begin{table}[htbp]
  \centering
  \small
  \caption{A worked example for three runs ($z^*=10$, $\tlim=30$). Lower is better. The DIMACS columns are shown on their conventional $\times100$ scale. Berthold's score is evaluated for Run~3 although the run passes below the reference at $t=8$, a case outside the measure's design envelope. A dash marks a value that does not exist: the raw gap has no finite whole-run average when the first incumbent arrives after $t=0$.}
  \label{tab:pi-example}
  \begin{tabular}{@{}lp{0.34\linewidth}rrrrrr@{}}
    \toprule
    & & & & & \multicolumn{2}{c}{DIMACS ($\times100$)} & \\
    \cmidrule(lr){6-7}
    Run & Valid incumbent events $(t_i,z_i)$ & $\ascore$ & $\ascore_B$ & Berthold & $\mathit{ACC}{=}1.1$ & $\mathit{ACC}{=}2$ & raw $\bar{g}$ \\
    \midrule
    1 & None & $1.00$ & $1.00$ & $1.00$ & $10.0$ & $100.0$ & --- \\
    2 & $(5,30),(7,29),(11,22),(14,18)$; \newline $(19,15),(24,14),(27,13)$ & $0.41$ & $0.53$ & $0.53$ & $10.0$ & $75.3$ & --- \\
    3 & $(0,14),(2,12),(3,11),(5,10),(8,9)$; \newline $(15,8),(18,7),(23,6),(26,5)$ & $-0.10$ & $-0.17$ & $0.23$ & $-18.3$ & $-16.0$ & $-0.16$ \\
    \bottomrule
  \end{tabular}
\end{table}

Run~2 exposes the effect of an acceptance threshold and of its placement. DIMACS at threshold 1.1 assigns it the same score as the empty run because no incumbent crosses $1.1z^*=11$, even though the valid objective improves by more than half. At threshold 2 the same rule admits the tail below $2z^*=20$ from $t=14$ onward and grades that progress 75.3 against the empty-run value of 100. Both signed threshold-free kernels grade the entire improvement.

The Berthold column evaluates his absolute gap inside its design envelope for Runs~1 and~2, where it coincides with $\ascore_B$. For Run~3 it is shown outside that envelope (degenerate under the analysis-time reference mode) where the absolute value folds the below-reference excursion back to a positive area: the resulting 0.23 no longer carries the sign of the improvement and reads like the score of a run held above the reference. This is no defect of Berthold's measure, whose analysis-time reference policy keeps the reference unbeaten by design. Considering the frozen-reference case is exactly what motivates one signed kernel covering the whole domain.

The raw-gap column makes the missing-run conventions concrete. A whole-run average of $g$ exists only for Run~3, which is feasible at the origin. Run~2 has no finite raw gap before its first incumbent and Run~1 has none at any time, so averaging them requires an external convention: dropping the run, assigning an arbitrary penalty, or flooring the start with a cheap initialization. The bounded kernels need no such choice, pre-incumbent states being simply their worst values.

Run~3 finally exposes the distinction between an endpoint and an area: its negative integrals record that the magnitude-weighted below-reference area dominates, while its raw final gap of $-50\%$ records the endpoint $z=5$. The printed cells also serve as fixed verification anchors for the released implementation.

\subsection{What aggregation preserves}

Scores are computed per run before any aggregation. For an expected-performance estimand, the arithmetic mean is natural and has the area identity used below. Signedness excludes ordinary geometric means across zero, but it does not make the arithmetic mean the only available summary: medians, quantiles, trimmed means, and other robust location estimators remain defined. They answer different questions and generally do not commute with integration.

The empirical campaign of Section~\ref{sec:evaluation} is organized in panels, each pairing an instance family with a size. Let $\mathcal{J}$ denote the set of panels, $\mathcal{I}_j$ the set of instances of panel $j\in\mathcal{J}$, and $\mathcal{R}_i$ the set of runs on instance $i\in\mathcal{I}_j$, with $\ascore_{jis}$ the score of the run of seed $s\in\mathcal{R}_i$ on instance $i$ of panel $j$. The empirical analysis uses
\begin{equation}
  \bar{\ascore}=\frac{1}{|\mathcal{J}|}\sum_{j\in\mathcal{J}}\frac{1}{|\mathcal{I}_j|}\sum_{i\in\mathcal{I}_j}\frac{1}{|\mathcal{R}_i|}\sum_{s\in\mathcal{R}_i}\ascore_{jis}.
  \label{eq:panel-equal}
\end{equation}
This panel-equal estimand gives each declared family-size panel equal headline weight. It is appropriate for the campaign design below, not universally preferable. A run-equal or instance-equal estimand may be correct for a different population. The weighting rule must be stated before interpreting any ranking.

\section{Reference Policies and Companion Interpretation}\label{sec:policies}\label{sec:kit}

\subsection{Analysis-time reference}\label{subsec:metric-classic}

In the classic regime, $z^*$ is the best value known after the experiment has completed. It may combine prior literature with the best solution found by any run. By construction, every scored incumbent satisfies $z(t)\ge z^*$, so both signed kernels and their integrals are nonnegative. Berthold's gap is used inside its intended design envelope, and $\sgap_B$ reproduces it exactly.

This policy answers a retrospective comparative question: given everything now known, how did each run progress toward the strongest available target? It has two benefits. The reference cannot be beaten within the analyzed data, simplifying interpretation, and all algorithms are normalized against the same strongest observed value. Its cost is temporal: a score cannot be finalized until the campaign and reference update are complete, and previously published scores change if they are recomputed against a stronger later reference. This mode is therefore appropriate for post-hoc analysis but it does not support online monitoring of a campaign and does not follow software release cycles or evolving versioned best-known-solution stores.

\subsection{Frozen, versioned reference}\label{subsec:metric-frozen}

In the frozen regime, $z^*$ comes from a named snapshot taken before the evaluated runs and remains immutable throughout execution and analysis. Modern software engineering practices have made versioning and reproducibility a first-class concern, and the same principles apply to best-known-solution stores. This is now an important feature of modern BKS stores. In the VRP literature, the MAMUT-routing repository \parencite{rascoussierMAMUTrouting2026}, the public reference repository accompanying \textcite{blauthVehicleRoutingTimedependent2024}, and the CombOpt historical records \parencite{comboptHistory} all provide versioned snapshots of best-known solutions, making new discoveries visible without silently altering prior scores. Modern tooling like Software Heritage \parencite{dicosmoSoftwareHeritageWhy2017} and Zenodo \parencite{zenodo2013} can archive and cite any snapshot of a BKS store, making it reconstructible rather than merely stable by convention. In this regime, the reference is not guaranteed to be unbeaten by the evaluated runs, producing negative instantaneous gaps when it happens. Freezing supports online campaign plot monitoring and comparison with previously computed scores. Versioning makes the anchor reconstructible rather than merely stable by convention.

The snapshot must carry per-instance values, provenance, a content digest, and the store version from which it was cut. Improvements are folded only in a separate post-processing step that creates a new snapshot. Scores computed against different snapshots are not silently mixed. This is especially important for best-known-solution repositories, where an entry is evidence of what was known at one date rather than a proof of optimality. Negative values also serve as useful diagnostics rather than certificates. They may indicate a real improvement, a stale reference, an instance mismatch, or a defect. A candidate apparently improving a certified optimum is an immediate integrity alarm. The sign preserves the event so that the analysis can investigate it, and does not decide which explanation is correct.

\begin{figure}[htbp]
  \centering
  \includegraphics[width=\linewidth]{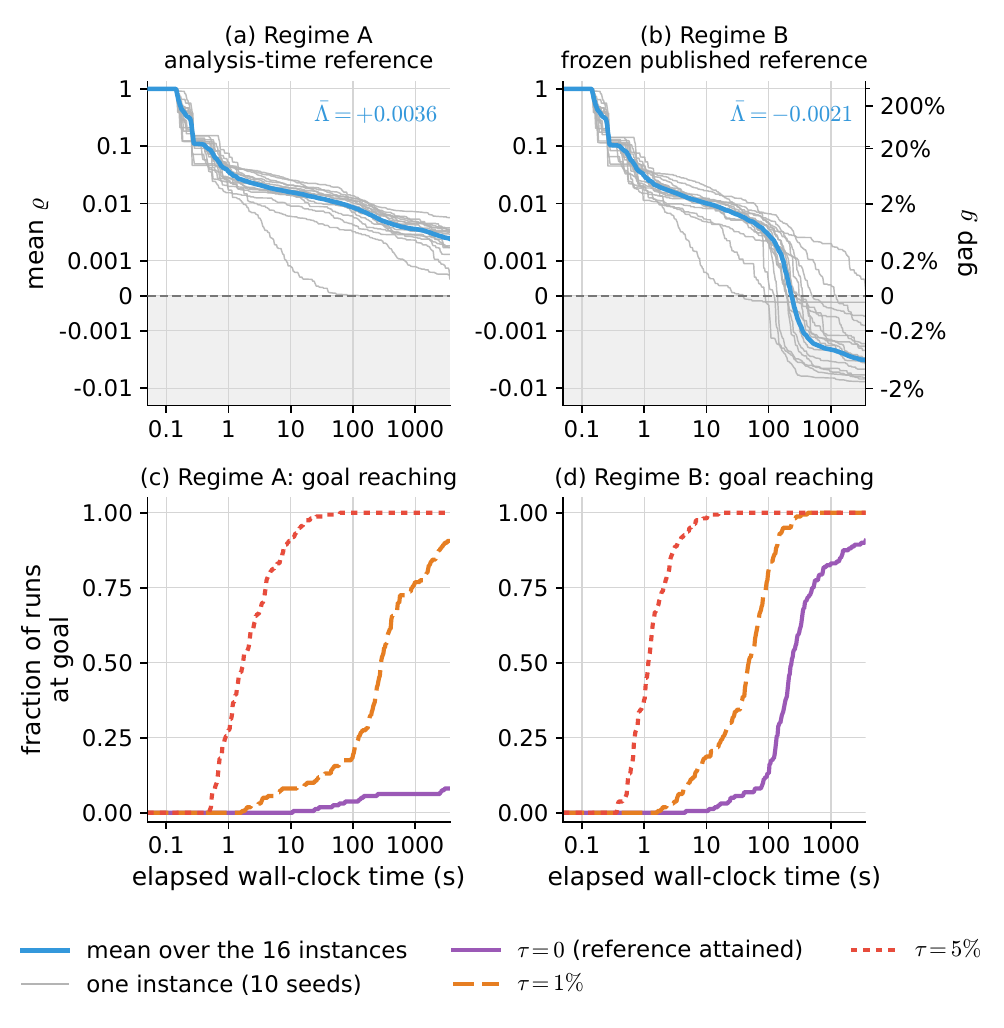}
  \caption{The same checker-valid trajectories under two reference policies. The analysis-time reference is the best value available after the campaign and keeps all gaps nonnegative. The frozen published snapshot may be beaten, so the mean curve and score can cross below zero. The panels isolate the effect of reference policy; they are not solver comparisons.}
  \label{fig:two-regime-campaign}
\end{figure}

Neither policy dominates. The analysis-time reference gives the cleanest retrospective comparison, while the frozen reference preserves historical comparability and makes discoveries visible. A paper should state which question it asks. Figure~\ref{fig:two-regime-campaign} scores the same trajectories twice and shows the semantic difference without changing a solver, instance, or timestamp.

\subsection{Recommended companion views}\label{subsec:kit-views}

\begin{figure}[htbp]
  \centering
  \includegraphics[width=\linewidth]{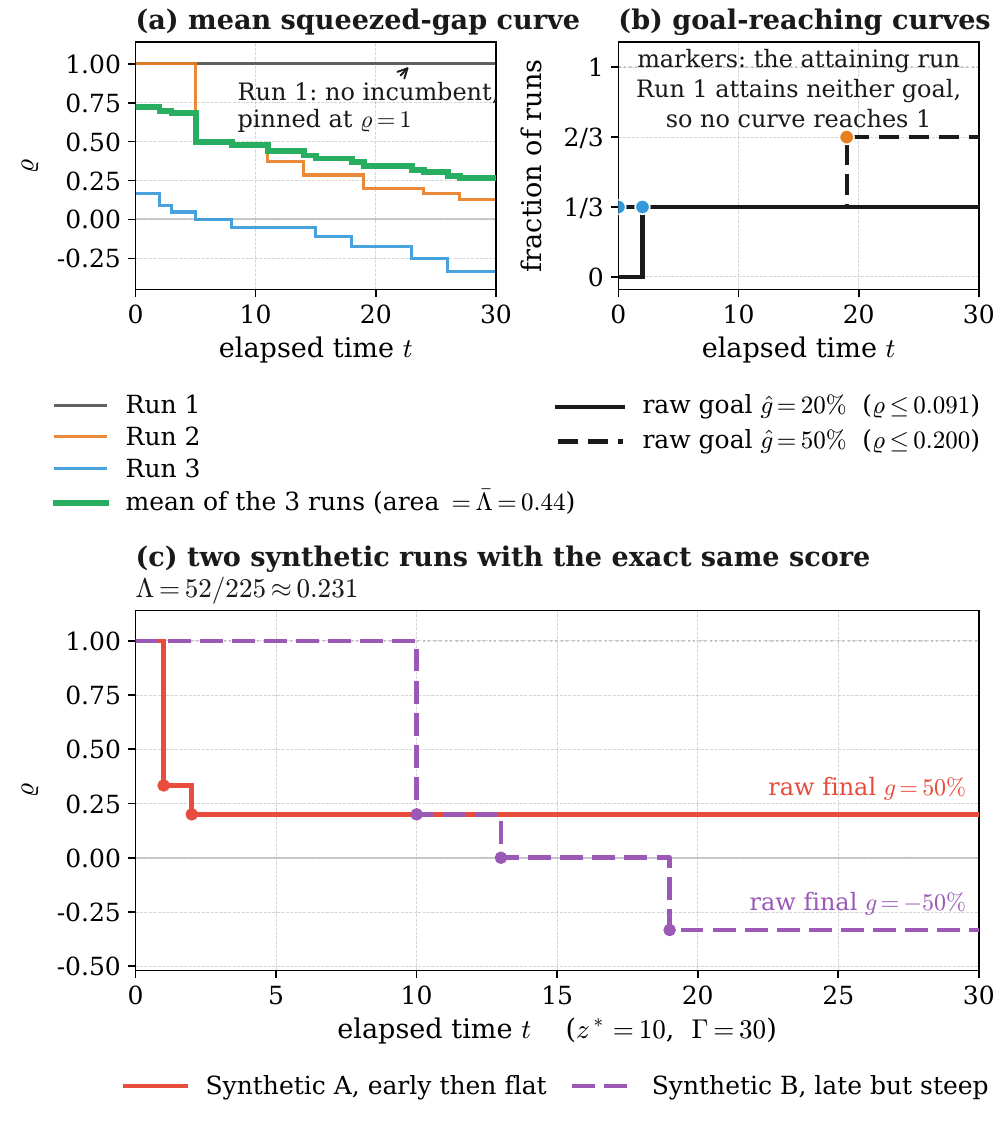}
  \caption{Why a whole-budget score benefits from companion views. Two trajectories can have equal normalized signed primal integrals while ending at different raw gaps and attaining declared targets at different times. The integral answers the area question, the raw final gap the endpoint question, and the attainment curve the rate-of-reaching question.}
  \label{fig:kit-companion}
\end{figure}

The scalar $\ascore$ describes whole-budget behavior but intentionally discards where improvements occurred. We therefore recommend reporting it with three familiar views whenever the study makes both whole-budget and endpoint claims: the raw final gap, a mean convergence curve, and one or more goal-attainment curves. Figure~\ref{fig:kit-companion} illustrates how equal integral scores can conceal different endpoints and attainment histories.

The three views are defined over a nonempty collection $\mathcal{R}$ of runs sharing the horizon $\tlim$; run $R\in\mathcal{R}$ carries its incumbent trajectory $z_R(t)$, its reference $z_R^*$, and hence its instantaneous squeezed gap $\sgap_R(t)=\sgap(z_R(t),z_R^*)$ and score $\ascore_R$. The reference is per-run because the collection may span several instances. The final gap is the familiar raw relative gap $g(z_R(\tlim),z_R^*)=(z_R(\tlim)-z_R^*)/z_R^*$; it states where a run ended and reveals a late improvement that may not change the integral's sign. A run still empty at $\tlim$ has no finite final gap, so its empty status and attainment contribution are reported rather than hidden inside a finite endpoint mean. Under a frozen reference, the signed final squeezed gap may additionally be useful when one bounded endpoint scale over all runs is desired. The mean convergence curve $\bar{\sgap}(t)=\frac{1}{|\mathcal{R}|}\sum_{R\in\mathcal{R}}\sgap_R(t)$ retains the temporal shape and remains finite before every run has an incumbent because empty states contribute 1. The goal-attainment curve for a raw gap goal $\hat{g}>0$ is
\begin{equation}
  A_{\hat{g}}(t)=\frac{1}{|\mathcal{R}|}\sum_{R\in\mathcal{R}}\mathbf{1}\!\left\{\min_{u\le t}\frac{z_R(u)-z_R^*}{z_R^*}\le\hat{g}\right\},
  \label{eq:attainment}
\end{equation}
where $\mathbf{1}\{\cdot\}$ is the indicator function and the minimum over an empty incumbent set is $+\infty$. An unattaining run contributes zero rather than disappearing. Goals are declared on the raw scale for readability and converted to their squeezed image $\hat{\sgap}=\hat{g}/(2+\hat{g})$ only inside the evaluator.

Linearity connects the score and mean curve. For the same collection of runs,
\begin{equation}
  \frac{1}{|\mathcal{R}|}\sum_{R\in\mathcal{R}}\ascore_R=\frac{1}{\tlim}\int_0^{\tlim}\frac{1}{|\mathcal{R}|}\sum_{R\in\mathcal{R}}\sgap_R(t)\,\mathrm{d}t=\frac{1}{\tlim}\int_0^{\tlim}\bar{\sgap}(t)\,\mathrm{d}t.
  \label{eq:area-identity}
\end{equation}
The identity provides an important implementation check. It holds because each run is represented on the same bounded scale and the arithmetic mean is linear. Confidence bands may be added to mean curves; an endpoint read from such a curve contains end-state information, but the explicit raw-final-gap distribution and attainment rate are usually easier to interpret.

\section{Empirical Evaluation}\label{sec:evaluation}\label{sec:applications}

\subsection{Research questions and data}

The empirical evaluation asks three questions. \textbf{RQ1, completeness and metric sensitivity:} how often do invalid or empty valid streams occur, and do the squeezed gap, signed max-form, Berthold gap, or DIMACS thresholds change scores and rankings? \textbf{RQ2, reference sensitivity:} what changes when identical trajectories are scored against a frozen pre-campaign snapshot or the post-campaign analysis reference? \textbf{RQ3, complementarity:} when do similar integral scores conceal different raw final gaps or target-attainment behavior?

The main corpus is a cross-solver campaign on the time-dependent vehicle routing problem with time windows (Duration-Minimization TDVRPTW) containing 8,800 one-hour runs: 212 instances in ten family-size panels, five solver arms, and five or ten seeds depending on instance size. The panels pair five instance families at two sizes each, from the classical time-dependent Solomon benchmark \parencite{dabiaBranchPriceTimeDependent2013} and the Lyon road-network instances of \textcite{rifkiImpactSpatiotemporalGranularity2020} to road-network families reaching one thousand customers. Four Tier-1 arms form the comparison set; a time-sliced Hexaly Tier-2 arm is descriptive because it changes the model representation. For transparency, the campaign originates in the evaluation of our own work: it measures the open-source KAYROS solver \parencite{rascoussierKAYROSTechReport2026} against established contenders from open-source and commercial organizations, Timefold \parencite{timefoldSolver2026}, jsprit \parencite{graphhopperJsprit2026}, and Hexaly \parencite{hexalyTDCVRPTW2026}. It thereby provides a concrete, complex environment in which to exercise the principles developed in this paper. Every captured incumbent was re-evaluated by an external checker on the original time-dependent functions. Two reference snapshots cover the same 212 instances: a pre-campaign snapshot frozen before any final run, and a post-fold analysis snapshot incorporating the best values found during the campaign. The instance set, checker, solver versions, and both snapshots are identified by content digests. The three research questions are descriptive post-campaign sensitivity analyses; the frozen confirmatory design concerned the original pooled solver contrasts, not the alternative rescoring or pair selection reported here. Table~\ref{tab:evaluation-design} summarizes the evidence used for each part of the analysis.

\begin{table}[htbp]
  \centering
  \small
  \caption{Evaluation design and the role of each artifact. Tier~1 supports equal-model solver comparisons. Tier~2 and the ladders diagnose the validity consequences of a time-sliced approximation.}
  \label{tab:evaluation-design}
  \begin{tabular}{@{}>{\raggedright\arraybackslash}p{0.22\linewidth}>{\raggedright\arraybackslash}p{0.21\linewidth}>{\raggedright\arraybackslash}p{0.22\linewidth}>{\raggedright\arraybackslash}p{0.27\linewidth}@{}}
    \toprule
    Evidence & Runs & Unit and horizon & Role in this paper \\
    \midrule
    Tier-1 campaign & $4\times1{,}760$ & 212 instances; one hour & Solver ordering, metric sensitivity, and reference sensitivity \\
    Tier-2 campaign & $1{,}760$ & Same instances and horizon & Descriptive failure-convention stress test \\
    Slice ladder & $3\times60$ & 6 instances, 10 seeds; one hour & Model-fidelity sensitivity across 5, 24, and 96 slices \\
    Thread ladder & $5\times100$ & 10 instances, 10 seeds; one hour & Validity under 1, 2, 4, 8, and 16 threads at 96 slices \\
    \bottomrule
  \end{tabular}
\end{table}

Two Hexaly Tier-2 ladders provide concentrated model-fidelity evidence \parencite{rascoussierHexalyThreadScaling2026}. Both use the natively time-sliced model, so neither depends on the Tier-1 evaluation path. The slice ladder contains 180 runs across three discretization levels. The selected 96-slice encoding is then evaluated in a 500-run thread ladder. Native feasibility in these models does not imply feasibility on the original functions, so the ladders expose the difference between captured candidates, checker-valid incumbent streams, invalid terminal solutions, and runs with no valid incumbent.

All 8,800 validated incumbent streams are re-analyzed by one deterministic procedure. It first reproduces the campaign's per-run $\ascore$ values and its pooled means, and only then rescores every trace under $\sgap$, $\sgap_B$, Berthold's absolute max-form, DIMACS thresholds 1.1 and 2, and both reference snapshots. Berthold's values under the frozen reference are out-of-envelope diagnostics, kept for completeness but not used in any ranking. Aggregation follows Equation~(\ref{eq:panel-equal}). The analysis reuses the frozen campaign records such that no new solver run was performed.

\subsection{RQ1: completeness and metric sensitivity}

Table~\ref{tab:health} shows why completeness is not a hypothetical concern. The four Tier-1 arms all end with valid terminal solutions and retain at least one checker-valid incumbent in every run. Hexaly Tier 1 nevertheless emits 109 captured incumbents rejected by the checker before its valid envelope is formed. Tier 2 emits 96,312 invalid captured incumbents. Of its runs, 79 end on an invalid terminal solution and 8 never produce any valid incumbent, and those 8 remain in every aggregate with gap 1 throughout.

\begin{table}[htbp]
  \centering
  \small
  \caption{Completeness of the 8,800-run campaign after independent checking. Each arm contains 1,760 runs. Invalid candidates are excluded from the best-so-far envelope, while a run with no valid incumbent scores 1 and remains in the denominator.}
  \label{tab:health}
  \begin{tabular}{@{}lrrr@{}}
    \toprule
    Arm & Invalid captured incumbents & Invalid terminal solutions & No valid incumbent \\
    \midrule
    KAYROS & 0 & 0 & 0 \\
    Timefold & 0 & 0 & 0 \\
    Hexaly Tier 1 & 109 & 0 & 0 \\
    jsprit & 0 & 0 & 0 \\
    Hexaly Tier 2 & 96,312 & 79 & 8 \\
    \bottomrule
  \end{tabular}
\end{table}

The ladders of Table~\ref{tab:ladder-validity} show how the failure pattern changes with discretization for the Hexaly Tier~2. Every level produces thousands of invalid captured incumbents, but only the two coarser settings leave runs with no valid incumbent at all, 12 of their 60 at five slices and 10 at twenty-four. The 96-slice setting leaves none, across both the slice and thread ladders. A bounded convention keeps these cells analyzable. The separate failure counts prevent that convenience from being mistaken for model fidelity.

The two signed kernels are empirically close in ordering and different in scale. The max-form roughly adds 50\% to each pooled mean (a factor of about 1.5) because its local scale is the raw gap rather than half of it (Table~\ref{tab:metric-pooled}; Appendix~\ref{app:empirical} tabulates both reference regimes), yet none of the 60 non-tied Tier-1 pairwise panel orderings reverses between the two kernels. At run level the integral sign can differ: under the frozen reference, 148 runs score negative under $\sgap$ against 166 under $\sgap_B$. These results support $\sgap$ as a coherent smooth working default on this campaign while leaving $\sgap_B$ as a valid alternative.

\begin{figure}[htbp]
  \centering
  \includegraphics[width=\linewidth]{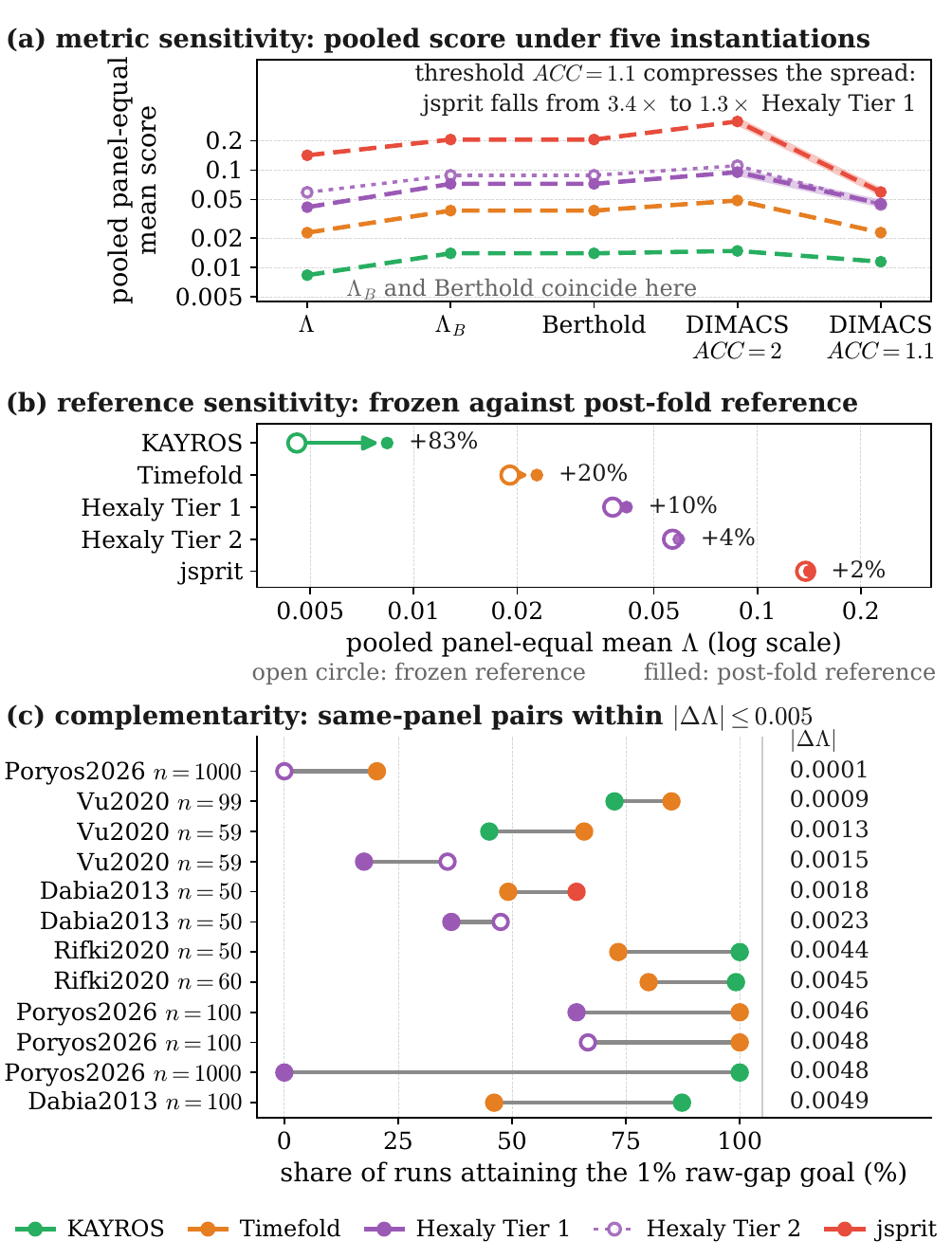}
  \caption{Sensitivity analysis over the complete 8,800-run campaign, one panel per research question, with Tier 2 descriptive throughout. (a) Pooled panel-equal scores under the five metric instantiations, post-fold reference: all agree on the ordering, but the DIMACS 1.1 threshold collapses the jsprit-to-Hexaly Tier 1 distance and reverses that pair on one panel. (b) Pooled $\ascore$ per arm under the frozen and post-fold snapshots: every score worsens against the stronger reference, and the ordering survives. (c) The 12 same-panel contrasts with similar whole-budget scores, the paired dots being the two arms' 1\% raw-gap attainment rates.}
  \label{fig:metric-sensitivity}
\end{figure}

The DIMACS acceptance threshold has a stronger consequence. Under the post-fold reference it assigns its flat worst score to 996 of 8,800 runs at threshold 1.1 and to 185 runs at threshold 2. The pooled Tier-1 ranking survives all five instantiations, but the threshold changes the distances that ranking is read from: jsprit stands 3.40 times above Hexaly Tier 1 under $\sgap$ and only 1.32 times above it at threshold 1.1, close enough that the two arms exchange places on the Lera2026 $n{=}200$ panel. That reversal is the only one the DIMACS instantiations produce among the 60 non-tied Tier-1 panel orderings, and it occurs under both reference snapshots. The one other discordance across the alternative instantiations is a near-tied pair that flips under the Berthold kernel in the frozen regime only (Appendix~\ref{app:empirical}). Figure~\ref{fig:metric-sensitivity}(a) traces the pooled scores across the five instantiations. The compression is not a defect in DIMACS, whose acceptance rule answers a competition-specific question. It demonstrates that the threshold is part of the estimand and cannot be treated as formatting.

\begin{table}[htbp]
  \centering
  \small
  \caption{Panel-equal pooled scores under the post-fold analysis reference. Lower is better. Here $\ascore_B$ equals Berthold's normalized primal integral because no incumbent beats the analysis reference. DIMACS values are divided by 100 to keep a common numerical scale. Tier~2 is descriptive.}
  \label{tab:metric-pooled}
  \begin{tabular}{@{}lrrrr@{}}
    \toprule
    Arm & $\ascore$ & $\ascore_B$ (Berthold) & DIMACS 1.1 & DIMACS 2 \\
    \midrule
    KAYROS & $0.00837$ & $0.01402$ & $0.01150$ & $0.01483$ \\
    Timefold & $0.02282$ & $0.03839$ & $0.02285$ & $0.04875$ \\
    Hexaly Tier 1 & $0.04162$ & $0.07221$ & $0.04526$ & $0.09507$ \\
    jsprit & $0.14160$ & $0.20521$ & $0.05953$ & $0.31459$ \\
    Hexaly Tier 2 & $0.05906$ & $0.08817$ & $0.04441$ & $0.11080$ \\
    \bottomrule
  \end{tabular}
\end{table}

\subsection{RQ2: reference sensitivity}

Folding the campaign's improvements into the reference store changes 54 of the 212 values. Rescoring identical trajectories against the stronger post-fold snapshot raises every pooled mean, because the same incumbents now sit farther from their anchors. Figure~\ref{fig:metric-sensitivity}(b) shows the shift for each arm, and Appendix~\ref{app:empirical} tabulates both regimes. The largest relative movement falls on the arms closest to the references, since a small absolute strengthening is a large share of a near-zero score.

The policy changes interpretation more than ordering in this dataset. The runs scoring negative under the frozen snapshot all become positive under the post-fold snapshot, because the improvements they represent have been folded into the reference. Yet none of the 60 non-tied Tier-1 pairwise panel orderings reverses between snapshots, and the pooled ordering is unchanged. This negative result is useful: reference provenance materially changes numerical values and whether discoveries remain visible, but the main comparison here remains robust.

The frozen reference preserves the historical question and exposes improvements made relative to that history. The analysis-time reference answers the cleaner retrospective question after folding. Publishing the snapshot digest makes both readings reproducible.

\subsection{RQ3: complementarity}

RQ3 asks what the companion views recover after the integral has summarized a trajectory. We screen the 100 within-panel arm pairs under the post-fold reference, including Tier~2 descriptively, and retain pairs with similar mean scores ($|\Delta\ascore|\le0.005$) but a difference of at least 0.01 in mean raw final gap or ten percentage points in attainment of the campaign's 1\% raw-gap goal. The goal was fixed before the campaign, but the pair-selection rule was developed afterward and is therefore exploratory. Twelve contrasts satisfy it. Appendix~\ref{app:rq3-sensitivity} lists them and shows that such cases remain under nearby thresholds, although their number changes. Figure~\ref{fig:metric-sensitivity}(c) displays the central observation directly. Each horizontal segment joins the attainment rates of two arms whose whole-budget scores are close, yet every segment spans at least ten percentage points. In the most pronounced case, KAYROS attains the 1\% goal in every run on Poryos2026 $n{=}1000$, whereas Hexaly Tier~1 never attains it, despite a difference of less than 0.005 between their mean integral scores.

The companion views confirm the interpretation of the score in this campaign: the arm with the better integral also attains the goal more often and, where comparable, finishes closer to the reference. The integral summarizes quality over the complete budget, while attainment preserves the probability of reaching a declared target. Their joint reading is therefore more informative than either alone, but the exploratory selection and the descriptive Tier~2 comparison limit this observation to the present campaign. Together, the results establish four campaign-level findings. (i) The signed-kernel choice preserves the observed solver ordering. (ii) An acceptance threshold compresses the distances between arms enough to reverse a panel ordering. (iii) Reference updates alter scale and sign but not that ordering. (iv) Companion views reveal distinctions the integral intentionally compresses. Such observations are specific to this campaign, but they illustrate the practical consequences of the evaluation choices. It shows that the signed integral is a coherent and informative summary of whole-budget behavior both in classic post-campaign and frozen pre-campaign reference regimes, while also demonstrating that it is best read with companion views that preserve endpoint and attainment information.

\section{Practical Use, Limitations, and Artifacts}\label{sec:practice}\label{sec:reproducibility}

\subsection{A compact reporting contract}

An anytime score is interpretable only together with the choices that produced it. A study should therefore state the objective, any positive transformation applied to it, and the independent validation procedure. It should declare the clock, the work included in the measured time, and the horizon. It should identify the reference policy, whether the analysis-time construction or a named frozen snapshot. It should give the gap kernel with its pre-incumbent value, and describe how invalid, crashed, and empty runs are treated and counted. A crash after valid events holds the last checker-valid incumbent to the horizon; a crash or missing output with no such event remains an empty stream, scores 1, and stays in every denominator. Finally, the study should state the aggregation order over seeds, instances, and panels, together with any attainment goals on the raw scale. Report $\ascore$ with raw final gap and curve views when both whole-budget and endpoint behavior matter.

Development and exploratory analysis may remain adaptive. When confirmatory statistical claims are made after choosing among many panels, budgets, metrics, and contrasts, freezing the final analysis plan before the final data helps distinguish confirmation from exploration. In the campaign reused here, instances are the bootstrap units, seeds remain attached to their instances, and Holm adjustment applies only to the three declared pooled contrasts \parencite{efronIntroductionBootstrap1993,holmSimpleSequentiallyRejective1979}. The alternative rescoring and RQ3 screen reported in this paper remain descriptive post-campaign analyses.

An arithmetic mean crossing zero should be reported through the signed difference in score units and, where helpful, the underlying component means. A relative percentage change is unstable or misleading when the baseline is near zero or when the sign changes. Medians and quantiles remain legitimate descriptive summaries, but the area identity in Equation~(\ref{eq:area-identity}) belongs specifically to the arithmetic expectation estimand.

\subsection{Limitations}

Strict positivity is a genuine representation constraint. Many maximization, zero-valued, or negative-valued objectives can be transformed into strictly positive minimization objectives, but additive shifts change the score and therefore become part of the evaluation contract. The transformation and its sensitivity must be declared with the results. A hierarchy or multiobjective problem likewise requires an explicit scalarization before any single integral can be computed. Appendix~\ref{app:lexico} records one calibrated hierarchical construction as an extended example; its component objectives still require separate reporting.

The score is scale-invariant but reference-sensitive, budget-conditional, and dependent on the recorded trace. Sparse polling can miss internal improvements, different initialization contracts can move early area, and a frozen best-known value may later be shown wrong. Protocol disclosure and a published reference snapshot delimit the resulting estimand. The empirical study is anchored in time-dependent vehicle routing. It covers only the solver families, budgets, objective, and checker represented there. The kernel and reference robustness results are descriptive of that domain. Different failure rates or quality distributions may produce ranking reversals between $\sgap$ and $\sgap_B$, which is why the implementation exposes both. The Tier-2 arm solves an approximation and is retained because it activates the failure conventions. It is not part of the equal-model confirmatory comparison. Counts of invalid candidates depend on capture frequency as well as model fidelity, so they should be read as evidence that invalidity is present, not as a normalized error rate between solvers.

\section{Conclusion}\label{sec:conclusion}

Anytime evaluation already has mature ingredients, but their free choices become consequential when feasibility is delayed, incumbents are invalid, or a frozen reference is beaten. The Normalized Signed Primal Integral studied here makes one coherent set of choices: a bounded threshold-free instantaneous gap, an intrinsic value for an empty stream, exact integration of the validated best-so-far trajectory, and a sign that remains meaningful under a frozen versioned reference. The signed max-form is a valid alternative, and our preference for the smooth difference-over-sum kernel is a default rather than a uniqueness claim.

The empirical results show both stability and sensitivity. The two signed kernels preserve all panel orderings in the considered experiment campaign, and strengthening references from the campaign results changes levels and signs without changing those orderings. A DIMACS acceptance threshold leaves the pooled ranking intact but shrinks the distance between two arms and reverses their order on one panel, while an exploratory screen identifies 12 panel contrasts in which similar integral scores conceal different endpoints or attainment rates. The practical lesson is to state the kernel and reference policy, retain every failed or empty run, and read the whole-budget score with the raw endpoint and target views needed by the claim.

Reference policy is the choice on which we take a position. A frozen snapshot precisely versions the selected per-instance reference set: each value can be traced to a named release or commit of an evolving BKS store such as MAMUT-routing, or to a pinned public benchmark repository. Newly discovered best solutions then produce an explicit new snapshot rather than silently changing the meaning of previously reported scores. The classic policy of folding the latest best-known values into the analysis remains the natural choice for a self-contained retrospective comparison, with a unified bounded metric working by design on both this classic setting as well as frozen versioned references. Implementation and artifacts are openly available and auditable.

\section*{Acknowledgements and Statements}

\paragraph{Acknowledgements.} We thank Adrien Pichon for the collaboration on MAMUT-routing and Romain Fontaine for his help with Grid'5000. Special thanks to Guillaume Beslon for his meaningful support. The experiments were carried out using the Grid'5000 testbed, supported by a scientific interest group hosted by Inria and including CNRS, RENATER, several universities, and other organizations (see \url{https://www.grid5000.fr}).

\paragraph{Funding.} This work was funded by the French National Research Agency (ANR) as part of the MAMUT project, ANR-22-CE22-0016, \enquote{Machine learning And Matheuristics algorithms for Urban Transportation}.

\paragraph{Data availability.} The estimator, the figure and analysis generators, the derived sensitivity data, the frozen instance set and reference snapshots, and the validated campaign streams form the replication package. The campaign archive is publicly deposited on Zenodo \parencite{rascoussierAnytimeSolverComparisonData2026}, and the analysis and plotting code that regenerates every table, statistic and figure from it is released at \url{https://github.com/0nyr/kayros-campaign-analysis}. The raw runs of the two Hexaly model-fidelity ladders of Section~\ref{sec:evaluation} are deposited separately as a public dataset \parencite{rascoussierHexalyThreadScalingData2026} accompanying version~2 of the thread-scaling report \parencite{rascoussierHexalyThreadScaling2026}. All conclusions and numerical results can be reproduced, audited or extended.

\paragraph{Use of generative AI.} Generative AI tools were used throughout the preparation of this work, including frontier models from Anthropic (Claude) and OpenAI (ChatGPT). They contributed to various aspects of this work, including code generation, experiment orchestration, initial drafting, and adversarial review of this paper. Florian Rascoussier led the review and revision of outputs from these tools. The authors take responsibility for the final content of the work.

\paragraph{Conflict of interest.} The authors declare no competing interests. One campaign compares a solver developed during the PhD research of Florian Rascoussier with independent contenders. Its instance set and confirmatory analysis were frozen before the final records existed, and every incumbent of every arm was re-evaluated by the same external checker.

%% file: sections/1-appendices.tex

\section{Mathematical Details}\label{app:math}

\subsection{The signed max-form variant}\label{app:maxform}

For $z,z^*>0$, let $x=z/z^*$. The signed max-form of Equation~(\ref{eq:maxform}) is
\begin{equation}
  \sgap_B(x)=
  \begin{cases}
    x-1, & 0<x<1,\\
    1-1/x, & x\ge1.
  \end{cases}
\end{equation}
It is strictly increasing, has range $(-1,1)$ on finite positive objectives, and tends to $1$ as $x\to+\infty$. The no-incumbent convention adds the endpoint 1. Reciprocal ratios are symmetric because $\sgap_B(1/x)=-\sgap_B(x)$. At $x=1$, the two first derivatives are both 1, while the second derivatives are 0 from the left and $-2$ from the right. The function is therefore $C^1$ with a curvature discontinuity.

Above the reference, $\sgap_B=1-z^*/z$, which is exactly Berthold's positive-objective gap in Equation~(\ref{eq:berthold-positive-gap}). Below the reference, $\sgap_B=(z-z^*)/z^*$, the raw excess gap. Its local expansion is $\sgap_B=g+\mathcal{O}(g^2)$ above the reference and exactly $g$ below it, whereas $\sgap=g/2+\mathcal{O}(g^2)$. The signed max-form consequently preserves the familiar local raw-gap scale.

The two kernels are related above the reference by
\begin{equation}
  \sgap_B=\frac{2\sgap}{1+\sgap}, \qquad \sgap=\frac{\sgap_B}{2-\sgap_B},
\end{equation}
and below it by
\begin{equation}
  \sgap_B=\frac{2\sgap}{1-\sgap}, \qquad \sgap=\frac{\sgap_B}{2+\sgap_B}.
\end{equation}
These monotone pointwise relations do not guarantee equal integral rankings because integration follows the nonlinear transformation. The empirical no-reversal result of Section~\ref{sec:evaluation} is therefore evidence rather than an algebraic necessity.

\subsection{Finite-sum evaluation and the area identity}

Let run $R\in\mathcal{R}$ have breakpoints $0=t_{R0}<t_{R1}<\cdots<t_{Rm_R}<t_{R,m_R+1}=\tlim$ after equal-time coalescing and best-so-far filtering, where $m_R$ is the number of retained incumbent events. Write $\sgap_{Ri}$ for the constant value the instantaneous squeezed gap $\sgap_R$ of Section~\ref{subsec:kit-views} holds on the half-open interval $[t_{Ri},t_{R,i+1})$, so that $\sgap_{R0}=1$ before the first valid incumbent. Equation~(\ref{eq:ascore}) is
\begin{equation}
  \ascore_R=\frac{1}{\tlim}\sum_{i=0}^{m_R}\sgap_{Ri}(t_{R,i+1}-t_{Ri}).
\end{equation}
It is exact because the recorded trajectory is a step function. To prove Equation~(\ref{eq:area-identity}), take the union of every run's breakpoints, refine each step function onto that shared partition, write $\Delta t_j$ for the width of its $j$-th cell and $\sgap_{Rj}$ for the value run $R$ holds there, and interchange the two finite sums:
\begin{align}
  \frac{1}{|\mathcal{R}|}\sum_{R\in\mathcal{R}}\ascore_R
  &=\frac{1}{|\mathcal{R}|\tlim}\sum_{R\in\mathcal{R}}\sum_j \sgap_{Rj}\Delta t_j\\
  &=\frac{1}{\tlim}\sum_j\left(\frac{1}{|\mathcal{R}|}\sum_{R\in\mathcal{R}}\sgap_{Rj}\right)\Delta t_j\\
  &=\frac{1}{\tlim}\int_0^{\tlim}\bar{\sgap}(t)\,\mathrm{d}t,
\end{align}
where $\bar{\sgap}$ is the mean convergence curve of Section~\ref{subsec:kit-views} over the runs of $\mathcal{R}$. The same proof holds for fixed nonnegative run weights summing to one. It does not extend to medians or quantiles because those operators are nonlinear.

\subsection{Raw and squeezed scales}

The inverse of Equation~(\ref{eq:squeezed-gap}) is $g=2\sgap/(1-\sgap)$. Table~\ref{tab:scale-landmarks} gives selected landmarks on the two scales.
\begin{table}[htbp]
  \centering
  \small
  \caption{Conversion between raw excess gaps and squeezed gaps.}
  \label{tab:scale-landmarks}
  \begin{tabular}{@{}rr@{}}
    \toprule
    Raw gap $g$ & Squeezed gap $\sgap$ \\
    \midrule
    $-50\%$ & $-33.33\%$ \\
    $-10\%$ & $-5.26\%$ \\
    $-5\%$ & $-2.56\%$ \\
    $0$ & $0$ \\
    $5\%$ & $2.44\%$ \\
    $10\%$ & $4.76\%$ \\
    $100\%$ & $33.33\%$ \\
    \bottomrule
  \end{tabular}
\end{table}

\section{Additional Empirical Results}\label{app:empirical}

\subsection{Pooled sensitivity table}

Table~\ref{tab:all-metrics} gives every pooled panel-equal mean of the sensitivity analysis. DIMACS values are expressed on the unitless gap scale rather than multiplied by 100. Berthold and $\sgap_B$ agree under the post-fold analysis reference because no incumbent lies below it. Under the frozen reference, the absolute max-form is evaluated only as an out-of-envelope diagnostic: it differs from $\sgap_B$ after an improvement and is not presented as Berthold's intended reference policy.

\begin{table}[htbp]
  \centering
  \small
  \caption{Pooled panel-equal means under the frozen pre-campaign reference (top) and post-fold analysis reference (bottom). Lower is better within a column; magnitudes are not compared across columns because the gap normalizations differ.}
  \label{tab:all-metrics}
  \begin{tabular*}{\linewidth}{@{\extracolsep{\fill}}lrrrrr@{}}
    \toprule
    Arm & $\sgap$ & $\sgap_B$ & Berthold & DIMACS 1.1 & DIMACS 2 \\
    \midrule
    \multicolumn{6}{c}{Frozen pre-campaign reference}\\
    KAYROS & $0.00457$ & $0.00675$ & $0.01085$ & $0.00462$ & $0.00702$ \\
    Timefold & $0.01905$ & $0.03154$ & $0.03770$ & $0.01748$ & $0.04045$ \\
    Hexaly Tier 1 & $0.03794$ & $0.06646$ & $0.06646$ & $0.04452$ & $0.08553$ \\
    jsprit & $0.13816$ & $0.20083$ & $0.20083$ & $0.05909$ & $0.30715$ \\
    Hexaly Tier 2 & $0.05664$ & $0.08450$ & $0.08450$ & $0.04379$ & $0.10537$ \\
    \midrule
    \multicolumn{6}{c}{Post-fold analysis reference}\\
    KAYROS & $0.00837$ & $0.01402$ & $0.01402$ & $0.01150$ & $0.01483$ \\
    Timefold & $0.02282$ & $0.03839$ & $0.03839$ & $0.02285$ & $0.04875$ \\
    Hexaly Tier 1 & $0.04162$ & $0.07221$ & $0.07221$ & $0.04526$ & $0.09507$ \\
    jsprit & $0.14160$ & $0.20521$ & $0.20521$ & $0.05953$ & $0.31459$ \\
    Hexaly Tier 2 & $0.05906$ & $0.08817$ & $0.08817$ & $0.04441$ & $0.11080$ \\
    \bottomrule
  \end{tabular*}
\end{table}

Under the post-fold reference, DIMACS 1.1 saturates at its worst value on 996 runs and DIMACS 2 on 185. Under the frozen reference the corresponding counts are 931 and 184. The 1.1 threshold preserves the pooled Tier-1 ranking in both regimes, but it brings jsprit within a factor 1.32 of Hexaly Tier 1 against a factor 3.40 under $\sgap$ on the post-fold reference (1.33 against 3.64 on the frozen one), and under both references it reverses the order of those two arms on the Lera2026 $n{=}200$ panel. That single reversal is the only discordance among the 60 non-tied Tier-1 panel orderings between $\sgap$ and any of the four alternative instantiations, apart from one further $\sgap$-against-Berthold reversal that arises only under the frozen reference, where Berthold's construction is applied outside its intended envelope.

\subsection{Hexaly Tier-2 ladder validity}

\begin{table}[htbp]
  \centering
  \small
  \caption{Checker validity in the Hexaly Tier-2 slice and thread ladders, both built on the natively time-sliced model. Each thread-ladder row contains 100 one-hour runs.}
  \label{tab:ladder-validity}
  \begin{tabular}{@{}lrrr@{}}
    \toprule
    Setting & Invalid captured incumbents & Invalid terminal solutions & No valid incumbent \\
    \midrule
    5 slices & 21,686 & 19 & 12 \\
    24 slices & 20,692 & 12 & 10 \\
    96 slices & 11,208 & 4 & 0 \\
    96 slices, 1 thread & 12,314 & 7 & 0 \\
    96 slices, 2 threads & 13,998 & 9 & 0 \\
    96 slices, 4 threads & 14,231 & 7 & 0 \\
    96 slices, 8 threads & 12,766 & 5 & 0 \\
    96 slices, 16 threads & 12,439 & 7 & 0 \\
    \bottomrule
  \end{tabular}
\end{table}

The invalid-candidate total is affected by the number of states captured and should not be interpreted as a normalized violation probability. The decisive completeness column is the number of runs with no checker-valid incumbent, which falls from 12 and 10 at the two coarser discretizations to zero at 96 slices.

\subsection{RQ3 threshold sensitivity}\label{app:rq3-sensitivity}

Table~\ref{tab:companion-contrasts} reports the 12 contrasts selected by the rule used in the core evaluation. In this appendix $g_{\mathrm{final}}$ denotes an arm's mean raw final gap $g(z_R(\tlim),z_R^*)$ over the runs of a panel, and $A_{1\%}$ its attainment rate $A_{\hat g}(\tlim)$ of Equation~(\ref{eq:attainment}) at the campaign goal $\hat g=0.01$, both as defined in Section~\ref{subsec:kit-views}.

\begin{table}[htbp]
  \centering
  \small
  \caption{All 12 within-panel contrasts selected by the exploratory complementarity rule. Differences are absolute, and attainment differences are percentage points.}
  \label{tab:companion-contrasts}
  \begin{tabular*}{\linewidth}{@{\extracolsep{\fill}}lllrrr@{}}
    \toprule
    Panel & Arm A & Arm B & $|\Delta\ascore|$ & $|\Delta g_{\mathrm{final}}|$ & $|\Delta A_{1\%}|$ \\
    \midrule
    Poryos2026 $n{=}1000$ & Timefold & Hexaly Tier 2 & $0.00012$ & $0.00328$ & $20.3$ \\
    Vu2020 $n{=}99$ & KAYROS & Timefold & $0.00089$ & $0.00278$ & $12.5$ \\
    Vu2020 $n{=}59$ & KAYROS & Timefold & $0.00135$ & $0.00299$ & $20.8$ \\
    Vu2020 $n{=}59$ & Hexaly Tier 1 & Hexaly Tier 2 & $0.00148$ & $0.00186$ & $18.3$ \\
    Dabia2013 $n{=}50$ & Timefold & jsprit & $0.00178$ & $0.00539$ & $15.0$ \\
    Dabia2013 $n{=}50$ & Hexaly Tier 1 & Hexaly Tier 2 & $0.00227$ & $0.00395$ & $10.8$ \\
    Rifki2020 $n{=}50$ & KAYROS & Timefold & $0.00440$ & $0.00608$ & $26.7$ \\
    Rifki2020 $n{=}60$ & KAYROS & Timefold & $0.00451$ & $0.00432$ & $19.2$ \\
    Poryos2026 $n{=}100$ & Timefold & Hexaly Tier 1 & $0.00460$ & $0.00722$ & $35.8$ \\
    Poryos2026 $n{=}100$ & Timefold & Hexaly Tier 2 & $0.00475$ & $0.00697$ & $33.3$ \\
    Poryos2026 $n{=}1000$ & KAYROS & Hexaly Tier 1 & $0.00484$ & $0.02361$ & $100.0$ \\
    Dabia2013 $n{=}100$ & KAYROS & Timefold & $0.00494$ & $0.00778$ & $41.3$ \\
    \bottomrule
  \end{tabular*}
\end{table}

Table~\ref{tab:rq3-threshold-sensitivity} varies each exploratory RQ3 threshold by a factor of two around the reported rule. The 1\% attainment goal itself was frozen before the campaign, but these pair-selection thresholds were chosen post campaign. The number of selected contrasts ranges from 2 under the tightest combination to 34 under the loosest. At the reported score-similarity cutoff $|\Delta\ascore|\le0.005$, varying the two material-difference cutoffs yields between 7 and 15 contrasts, so the existence of complementarity cases does not depend on the exact central values even though their count does.

\begin{table}[htbp]
  \centering
  \small
  \caption{Number of RQ3 contrasts selected among 100 within-panel arm pairs across a factor-two threshold grid. Columns vary the minimum difference in 1\% raw-gap attainment; rows vary the maximum score difference and minimum raw-final-gap difference. The bold entry is the rule reported in the core.}
  \label{tab:rq3-threshold-sensitivity}
  \begin{tabular}{@{}ccrrr@{}}
    \toprule
    Maximum $|\Delta\ascore|$ & Minimum $|\Delta g_{\mathrm{final}}|$ & $5$ pp & $10$ pp & $20$ pp \\
    \midrule
    $0.0025$ & $0.005$ & 8 & 7 & 4 \\
    $0.0025$ & $0.010$ & 8 & 6 & 2 \\
    $0.0025$ & $0.020$ & 8 & 6 & 2 \\
    \midrule
    $0.0050$ & $0.005$ & 15 & 13 & 9 \\
    $0.0050$ & $0.010$ & 15 & \textbf{12} & 7 \\
    $0.0050$ & $0.020$ & 15 & 12 & 7 \\
    \midrule
    $0.0100$ & $0.005$ & 34 & 32 & 26 \\
    $0.0100$ & $0.010$ & 34 & 31 & 24 \\
    $0.0100$ & $0.020$ & 33 & 30 & 23 \\
    \bottomrule
  \end{tabular}
\end{table}

\section{Logarithmic Time Reweightings}\label{app:logtime}

Equation~(\ref{eq:ascore}) averages uniformly in wall time. Equivalently, it reports the expected instantaneous gap when the interruption time is uniformly distributed over the declared budget. This is a natural default when each unit of elapsed time has the same operational importance, but it is not the only possible reading. Some applications care primarily about how quickly a solver becomes useful, while others emphasize the quality available immediately before a hard deadline. Logarithmic time scales provide two corresponding views and are familiar in anytime and black-box optimization evaluation \parencite{bertholdMeasuringImpactPrimal2013,bertholdConfinedPrimalIntegral2021,hansenCOCOPlatformComparing2021}.

A general weighted member of the same family is
\begin{equation}
  \ascore_w=\frac{\int_0^\tlim w(t)\sgap(z(t),z^*)\,\mathrm{d}t}{\int_0^\tlim w(t)\,\mathrm{d}t},
  \label{eq:ascore-weighted}
\end{equation}
for a nonnegative integrable weight $w$. The denominator normalizes the expression, so $\ascore_w$ remains a weighted average of instantaneous squeezed gaps and keeps the range $(-1,1]$. The uniform weight $w\equiv1$ gives $\int_0^\tlim w=\tlim$ and returns Equation~(\ref{eq:ascore}) exactly, so the score studied in this paper is the uniform member of this family rather than an object outside it.

Normalizing the weight into the probability density $f=w/\int_0^\tlim w$ on $[0,\tlim]$ gives the reading that makes the family interpretable,
\begin{equation}
  \ascore_w=\mathbb{E}_{T\sim f}\left[\sgap(z(T),z^*)\right],
  \label{eq:ascore-interruption}
\end{equation}
the expected instantaneous gap when the run is interrupted at a random time $T$ drawn from $f$. The weight is therefore a model of when the budget is cut short, and the uniform choice states that every moment of the declared budget is equally likely to be the one that matters.

Three consequences organize the rest of this appendix. First, the two logarithmic measures below are choices of $w$ inside one family rather than separate constructions, and they are read against $\ascore$ on that basis. Second, because $w$ is the interruption-time distribution, replacing it poses a different operational question instead of correcting a bias in the uniform score, so neither logarithmic variant supersedes $\ascore$. Third, the cutoffs both variants require are mandatory rather than cosmetic: $1/t$ is not integrable at the origin and $1/(\tlim-t)$ is not integrable at the deadline, so the normalizer of Equation~(\ref{eq:ascore-weighted}) diverges and the domain must be truncated. The truncation point defines $f$ and therefore belongs to the estimand, which is why it must be declared with any score computed under such a weight.

An early-weighted logarithmic average uses $w(t)=1/t$ on $[t_{\min},\tlim]$, for a cutoff $0<t_{\min}<\tlim$:
\begin{equation}
  \ascore_{\log}=\frac{1}{\ln(\tlim/t_{\min})}\int_{t_{\min}}^\tlim\frac{\sgap(z(t),z^*)}{t}\,\mathrm{d}t.
  \label{eq:ascore-log}
\end{equation}
This measure assigns equal importance to equal multiplicative intervals of elapsed time, so each order of magnitude contributes equally. It is useful when startup responsiveness matters, for example when solutions are consumed interactively, runs are frequently interrupted, or budgets span several orders of magnitude. Its cutoff is the first instance of the third consequence above: moving $t_{\min}$ toward zero gives increasing weight to initialization and to the period before the first valid incumbent, so the reported score depends on a value the timing contract must fix.

\begin{figure}[htbp]
  \centering
  \includegraphics[width=\linewidth]{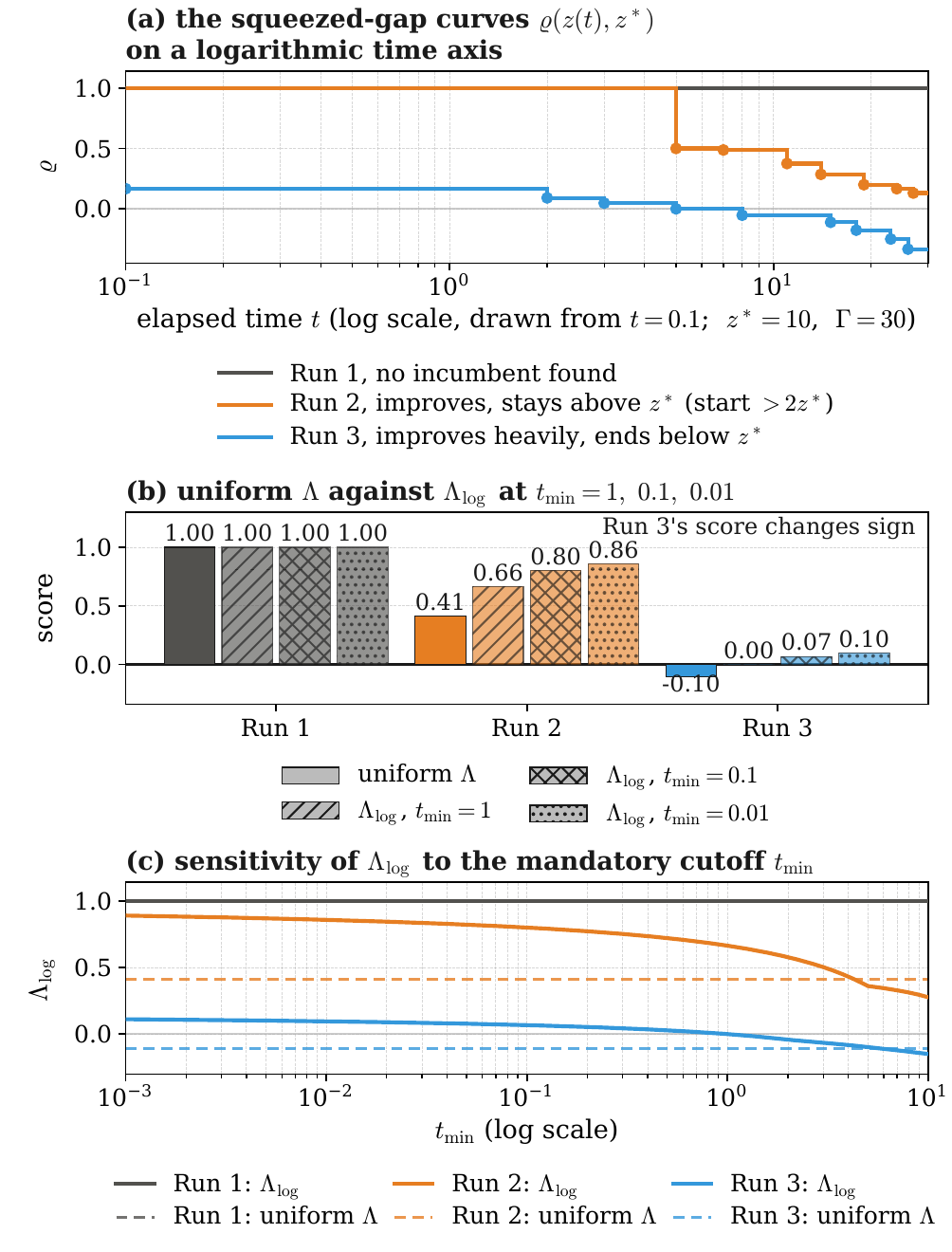}
  \caption{The early-weighted logarithmic average on the worked-example runs ($z^*=10$, $\tlim=30$). (a) The squeezed-gap trajectories on a logarithmic time axis. (b) The uniform $\ascore$ against $\ascore_{\log}$ at three cutoffs: the cutoff shifts every score and changes the sign of Run 3's. (c) Sensitivity of $\ascore_{\log}$ to the mandatory cutoff $t_{\min}$, with the cutoff-free uniform scores as dashed baselines.}
  \label{fig:logtime-early}
\end{figure}

Figure~\ref{fig:logtime-early} illustrates both the intended effect and its cost. The slow start of Run~2 becomes much more prominent, which is precisely the distinction sought from an early-weighted view. For Run~3, the same emphasis gives its brief initial period above the reference more influence than its later improvement below the reference, changing the sign of the integral. Panel~(c) further shows that the conclusion depends on the chosen cutoff. Early logarithmic weighting can therefore be informative as a diagnostic of responsiveness, but its score is no longer directly comparable with a uniform-time score. When only a visual account of early behavior is needed, plotting the ordinary convergence curve on a logarithmic time axis exposes the same temporal structure without changing the scalar estimand.

The end-weighted mirror applies the logarithm to the remaining time $\tlim-t$ and uses a cutoff $t_{\mathrm{cut}}\in(0,\tlim)$ on that remaining time, the mirror of $t_{\min}$:
\begin{equation}
  \ascore_{\mathrm{end}}=\frac{1}{\ln(\tlim/t_{\mathrm{cut}})}\int_0^{\tlim-t_{\mathrm{cut}}}\frac{\sgap(z(t),z^*)}{\tlim-t}\,\mathrm{d}t.
  \label{eq:ascore-end}
\end{equation}
Equal multiplicative intervals of remaining time now receive equal importance, so the final part of the run is expanded rather than the beginning. This view can be useful when a fixed deadline is operationally dominant or when the question concerns a solver's ability to exploit its last portion of budget. As $t_{\mathrm{cut}}\to0$, the score converges to the final squeezed gap whenever the trajectory is constant near the deadline. The measure is thus a parameterized blend of whole-run and endpoint performance rather than a distinct replacement for either.

\begin{figure}[H]
  \centering
  \includegraphics[width=\linewidth]{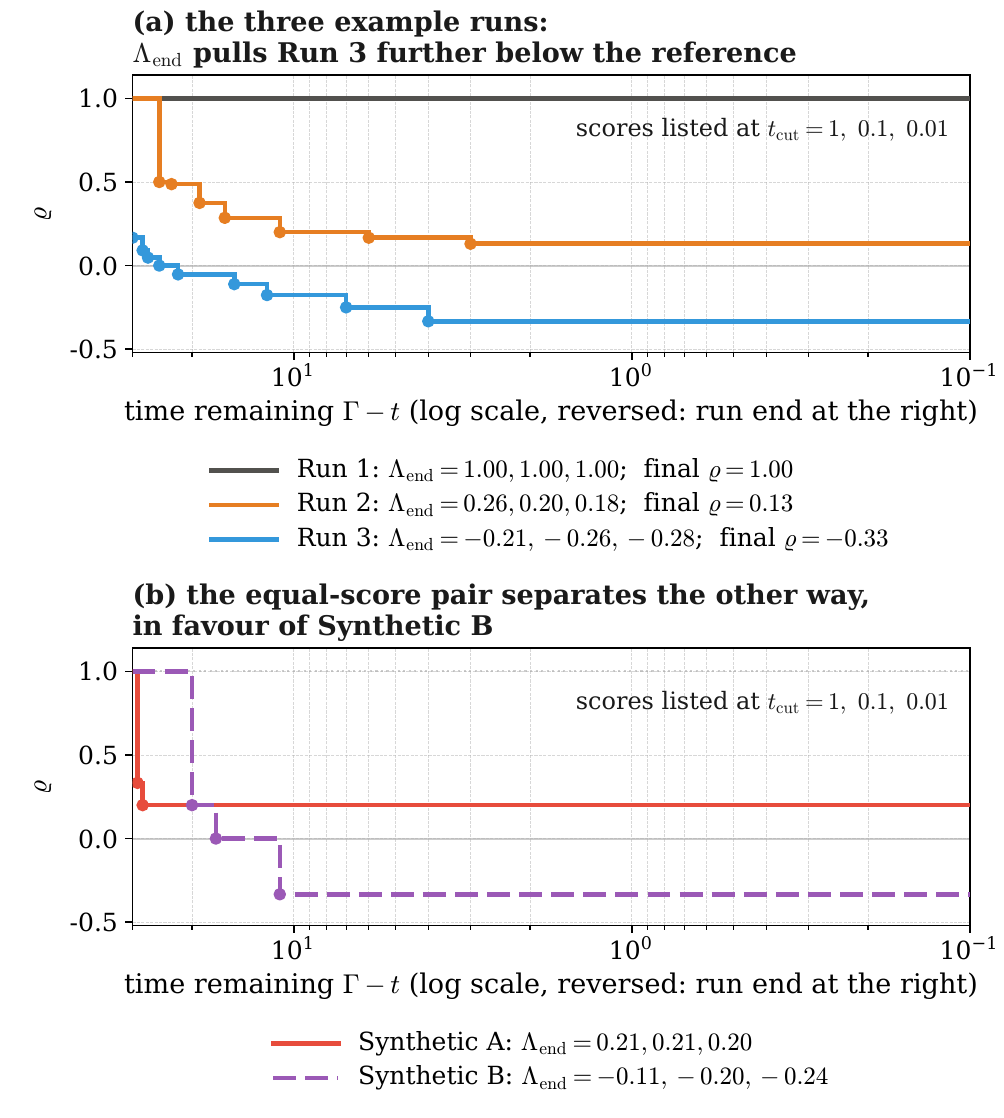}
  \caption{The end-weighted mirror on a reversed logarithmic remaining-time axis. (a) The worked-example runs: $\ascore_{\mathrm{end}}$ pulls Run 3 toward its final squeezed gap as $t_{\mathrm{cut}}$ shrinks. (b) A synthetic pair with identical uniform scores separates under the end-weighted view, in favor of the late, steep improver.}
  \label{fig:logtime-end}
\end{figure}

Figure~\ref{fig:logtime-end} makes this interpolation visible. Run~3 is pulled toward its better final gap as the cutoff shrinks. The synthetic pair in panel~(b) has identical uniform integrals, yet the end-weighted view favors the trajectory that improves late and finishes better. This distinction may be relevant under a hard deadline, but reporting the uniform integral and final gap separately is usually easier to interpret because it reveals whether whole-budget or endpoint behavior drives the comparison.

The two reweightings therefore have clear but limited roles. Early weighting asks which solver becomes useful sooner; end weighting asks which solver makes best use of the deadline. Either can serve as a secondary scalar when that question is explicit, provided its cutoff is declared as part of the timing contract. For general reporting, the uniform $\ascore$ remains preferable because it introduces no cutoff and retains its average-over-the-budget interpretation. Early or late emphasis can then be supplied transparently through convergence curves on an appropriate time axis and through the separately reported final gap.

\section{Hierarchical Objectives as an Extended Example}\label{app:lexico}

Some vehicle-routing benchmarks minimize a primary integer resource count (typically the number of vehicles) and then a secondary continuous cost (like the total travel time or distance), written $\left| \mathcal{S} \right|$ and $C(\mathcal{S})$ for a feasible solution $\mathcal{S}$, the fleet size being the cardinality of the set of routes rather than a symbol of its own. A single trajectory score requires a scalar objective, commonly
\begin{equation}
  z(\mathcal{S})=M\left| \mathcal{S} \right|+C(\mathcal{S}),
  \label{eq:lexico-objective}
\end{equation}
where the scalarization multiplier $M>0$ prices one unit of the primary resource in units of the secondary cost. The multiplier must belong to the evaluation contract rather than to a solver's internal search penalty. Otherwise solvers are evaluated on different objective scales.

For two solutions with $\left| \mathcal{S}' \right|=\left| \mathcal{S} \right|+1$ (one more vehicle), Equation~(\ref{eq:lexico-objective}) preserves the lexicographic preference for $\mathcal{S}$ precisely when $M>C(\mathcal{S})-C(\mathcal{S}')$. A sufficient condition over an evaluated instance set is therefore
\begin{equation}
  M>\max\left\{C(\mathcal{S})-C(\mathcal{S}'):\left| \mathcal{S}' \right|=\left| \mathcal{S} \right|+1\right\}.
  \label{eq:lexico-condition}
\end{equation}
The condition must be checked on the relevant family. The convention $M=10n$, with $n$ the number of customers of the instance is an empirical calibration, not a universal constant.

A large multiplier creates a second interpretability problem even when it preserves order: changes in the secondary cost become nearly invisible beside the primary term. The normalized signed primal integral on $z$ must therefore be reported with the two objective components separately. This appendix illustrates how a scalar metric can be extended after an evaluation objective has been declared; it is not a contribution to hierarchical scalarization.